\documentclass{amsart}
\usepackage[dvips]{graphicx}
\usepackage{enumerate}
\usepackage{amscd}
\usepackage{amssymb}

\title{Geometric limits of quasi-Fuchsian groups}
\author{Teruhiko Soma}
\thanks{This paper is a simplified and cleared version of already circulated ones.
Any comments or questions on this paper are welcome.}

\address{Department of Mathematics and Information Sciences,
Tokyo Metropolitan University,
Minami-Ohsawa 1-1, Hachioji, Tokyo 192-0397, Japan}

\email{tsoma@center.tmu.ac.jp}

\begin{document}

\maketitle

\begin{abstract}
In this paper, we will determine the topological types of hyperbolic 3-manifolds $\mathbf{H}^3/G$ such that $G$ is a geometric limit of any algebraically convergent sequence of quasi-Fuchsian groups.
\end{abstract}

\newtheorem{theorem}{Theorem}[section]
\newtheorem{cor}[theorem]{Corollary}
\newtheorem{lemma}[theorem]{Lemma}
\newtheorem{sublemma}[theorem]{Sublemma}
\newtheorem{prop}[theorem]{Proposition}

\newtheorem{mtheorem}{Theorem}
\renewcommand{\themtheorem}{\Alph{mtheorem}}
\newtheorem{mcorollary}[mtheorem]{Corollary}

\theoremstyle{definition}
\newtheorem{definition}[theorem]{Definition}
\newtheorem{example}[theorem]{Example}
\newtheorem{xca}[theorem]{Exercise}
\newtheorem{remark}[theorem]{Remark}
\newtheorem{problem}[theorem]{Problem}

\numberwithin{figure}{section}
\numberwithin{equation}{section}

\def\nn{\mathbf{N}}
\def\rr{\mathbf{R}}
\def\cc{\mathbf{C}}

\def\cb{\mathcal{B}}
\def\ce{\mathcal{E}}
\def\ck{\mathcal{K}}
\def\cf{\mathcal{F}}
\def\cg{\mathcal{G}}
\def\ch{\mathcal{H}}
\def\cp{\mathcal{P}}
\def\cx{\mathcal{X}}
\def\cn{\mathcal{N}}
\def\cs{\mathcal{S}}
\def\cu{\mathcal{U}}
\def\cy{\mathcal{Y}}
\def\cz{\mathcal{Z}}
\def\eset{\emptyset}
\def\part{\partial}
\def\fd{\pi_1}
\def\ol{\overline}
\def\Int{\mathrm{Int}}
\def\fr{\mathrm{Fr}}
\def\thin{{\mathrm{thin}}}
\def\thick{{\mathrm{thick}}}
\def\pthin{{\mathrm{p:thin}}}
\def\pthick{{\mathrm{p:thick}}}
\def\ppthin{{\mathrm{(p):thin}}}
\def\ppthick{{\mathrm{(p):thick}}}
\def\to{\longrightarrow}
\def\sto{\rightarrow}
\def\ve{\varepsilon}
\def\L{\Lambda}
\def\O{\Omega}
\def\G{\Gamma}
\def\D{\Delta}
\def\Sg{\Sigma}

Let $\rho_n:\Pi\to \mathrm{PSL}_2(\cc)$ $(n\in \nn)$ be Kleinian representations 
from a torsion-free non-abelian group $\Pi$.
If $\{\rho_n\}$ converges algebraically to a representation $\rho_\infty$, then the sequence 
$\{\G_n\}$ with 
$\G_n=\rho_n(\Pi)$ has a subsequence converging geometrically to a Kleinian group $G$ containing $\G_\infty =
\rho_\infty(\Pi)$.
Various results on the relationship between $\G_\infty$ and $G$ have been obtained by Thurston \cite[Chapters 8, 9]{tl}, Kerckhoff and Thurston \cite{kt}, Ohshika \cite{oh1,oh2}, Anderson and Canary \cite{ac1,ac2,ac3}, Brock \cite{br} and so on.
Here, we are especially interested in geometric limits of quasi-Fuchsian groups, which are fundamental geometrically finite Kleinian groups.
Even if two sequences $\{\rho_n\}$, $\{\rho_n'\}$ of quasi-Fuchsian representations have the 
same algebraic limit, the sequences $\{\rho_n (\Pi)\}$, $\{\rho_n'(\Pi)\}$ do not necessarily have the same geometric limit.
This suggests that geometric limits have data which are not obtained from any data of algebraic limits.
In order to proceed furthermore with the study of geometric limits $G$ of quasi-Fuchsian groups, it would be meaningful for us 
to understand the topological types of the hyperbolic 3-manifolds $\mathbf{H}^3/G$.

Let $\Sg$ be a closed orientable surface of genus greater than one.
We fix a hyperbolic structure on $\Sg$ for convenience, and set $\Pi=\fd(\Sg)$.
In \cite{jm}, J{\o}rgensen and Marden gave an example of faithful representations $\zeta_n:\mathbf{Z}\to \mathrm{PSL}_2(\cc)$ with $\zeta_n(1)$ loxodromic such that the cyclic Kleinian groups $\zeta_n(\mathbf{Z})$ converge geometrically to a rank two parabolic group.
This is one of typical phenomena which appear in geometric limits.
In fact, Kerckhoff and Thurston \cite{kt} considered the cyclic action on the Bers slice $B_{\sigma_+}=\mathrm{Teich} (\Sg)\times \{\sigma_+\}$ at $\sigma_+\in \mathrm{Teich} (\Sg)$ generated by the Dehn twist $\varphi$ on $\Sg$ along a simple closed geodesic $l$.
Then they showed that any geometric accumulation point of the orbit $\{(\varphi_*^n(\sigma_-),\sigma_+)\}\subset B_{\sigma_+}$ is a Kleinian group $G$ such that $\mathbf{H}^3/G$ is homeomorphic to $\Sg\times (0,1)- l\times \{1/2\}$ for any $\sigma_-\in \mathrm{Teich}(\Sg)$.
Here a tubular neighborhood of $l\times \{1/2\}$ in $\Sg\times (0,1)$ corresponds to a $\mathbf{Z}\times \mathbf{Z}$-cusp of $\mathbf{H}^3/G$ where 
J{\o}rgensen-Marden phenomenon occurs.
By using this method iteratively, it is also possible to construct an example of a geometric limit $G'$ of quasi-Fuchsian groups such that $\mathbf{H}^3/G'$ has infinitely many $\mathbf{Z}\times \mathbf{Z}$-cusps.
In particular, $G'$ is infinitely generated.
Another important example of geometric limits of quasi-Fuchsian groups is given by Brock \cite{br}.
He considered the cyclic action on a Bers slice generated by a homeomorphism $\psi:\Sg\to \Sg$ such that $\psi|\Int H:\Int H\to \Int H$ is pseudo-Anosov and $\psi|(\Sg-\Int H)$ is the identity for a proper subsurface $H$ of $\Sg$.
Then any geometric accumulation point of the cyclic orbit $\{(\psi_*^n(\sigma_-),\sigma_+)\}\subset B_{\sigma_+}$ is a Kleinian group $G''$ such that $\mathbf{H}^3/G''$ is homeomorphic to $\Sg\times (0,1)-H\times \{1/2\}$.
However, all of these examples are very special ones.
In this paper, we will present what kinds of topological types appear generally in geometric limits of quasi-Fuchsian groups.

Let $\mathrm{pr}_\mathrm{hz}:\Sg\times I\to \Sg$ and $\mathrm{pr}_\mathrm{vt}:\Sg\times I\to I$ be the projections onto the first and second factors, where $I$ is the closed interval $[0,1]$.
For any $y\in I$, the preimage $\Sg_y=\mathrm{pr}_\mathrm{vt}^{-1}(y)$ is supposed to have the hyperbolic structure such that 
the restriction $\mathrm{pr}_\mathrm{hz}|\Sg_y:\Sg_y\to \Sg$ is an isometry.

Let $\cx$ be the closed subset of $\Sg\times I$ given below.
Set  $\cy=\mathrm{pr}_\mathrm{vt}(\cx)$, $X_y=\Sg_y\cap \cx$ for $y\in \cy$, $\L_{y}^+=\Sg_y\cap \ol{(\Sg\times (y,1]\cap \cx)}$ for $y<1$ and $\L_{y}^-=\Sg_y\cap \ol{(\Sg\times [0, y)\cap \cx)}$ for $y>0$.
Consider a maximal open geodesic subsurface $Z(\L_y^\ve)$ of $X_y$ disjoint from $\L_y^\ve$, which is determined uniquely by Lemma \ref{l22}.
Let $\lambda(\L_y^\ve)$ be the union of all simple geodesic loops $l$ in $X_y\setminus Z(\L_y^\ve)$ each of which 
has a \emph{free side} with respect to $\L_y^\ve$, that is, for the $\delta$-neighborhood 
$\cn_\delta(l,\Sg_y)$ with sufficiently small $\delta >0$, at least one component of 
$\cn_\delta (l,\Sg_y)\setminus l$ is disjoint from $\L_y^\ve$, see Fig.\ \ref{d_l} in Subsection \ref{p_c}.
From the definition, we know that $\lambda(\L_y^\ve)$ contains the frontier $\fr(X_y)$ of $X_y$ in $\Sg_y$.
Moreover, by Lemma \ref{l23}, $\lambda(\L_y^\ve)$ is a disjoint union of geodesic loops.

\begin{mtheorem}\label{main_a}
Let $\{\rho_n:\Pi\to \mathrm{PSL}_2(\cc)\}_{n=1}^\infty$ be any algebraically convergent sequence of quasi-Fuchsian representations such that $\{\rho_n(\Pi)\}_{n=1}^\infty$ converges geometrically to a Kleinian group $G$.
Then the Kleinian $3$-manifold $\mathbf{H}^3\cup \O(G)/G$ is homeomorphic to $\Sg\times I\setminus \cx$ such that $\cx$ is a closed subset of $\Sg\times I$ satisfying the following conditions {\rm (i)}-{\rm (iii)}.
\begin{enumerate}[\rm (i)]
\item
$\Sg\times I\setminus \cx$ is connected and contains $\Sg_{1/2}$.
\item
For any $y\in \cy$, $X_y$ is a disjoint union of a geodesic subsurface and simple geodesic loops in $\Sg_y$.
Moreover, for $\varepsilon =\pm$, each component of $Z(\L_{y}^{\varepsilon})$ is not homeomorphic to an open pair of pants.
\item
For any $y,z\in \cy$ with $y<z$, if a component $l_y$ of $\lambda(\L_{y}^+)$ is parallel to a component $l_z$ of 
$\lambda(\L_{z}^-)$ in $\Sg\times I\setminus \cx$, then $l_y$ and $l_z$ are horizontally parallel to each other in $\cx$.
\end{enumerate}
\end{mtheorem}

Theorem \ref{main_a} is proved by using the fact that $\mathbf{H}^3\cup \O(G)/G$ admits a block decomposition (see Proposition 
\ref{M_block}).
However, we should remark that an open set in $\Sg\times I$ admitting a block decomposition does not necessarily 
have a hyperbolic structure.
In fact, if a closed subset of $\Sg\times I$ satisfies the conditions (i) and (ii) but 
does not (iii), then the complement can not be a hyperbolic manifold. 

Any component of $Z(\L_{y}^{\varepsilon})$ (resp.\ of $\Sg_0\cup \Sg_1\setminus \cx$) 
is called a \emph{non-peripheral} (resp.\ \emph{peripheral}) \emph{front} of $\cx$.
It follows from Theorem \ref{main_a} that $\O(G)/G$ is homeomorphic to the union $\Sg_0\cup \Sg_1\setminus \cx$ of peripheral fronts.
Though in general the hyperbolic 3-manifold $\mathbf{H}^3/G$ has infinitely many ends, this fact shows that the end of $\mathbf{H}^3/G$ 
corresponding to any non-peripheral front of $\cx$ is geometrically infinite.

From now on, we call that any closed subset $\cx$ of $\Sg\times I$ satisfying 
the conditions (i)-(iii) in Theorem \ref{main_a} is \emph{crevasse} of $\Sg\times I$.

\begin{remark}
In particular, Theorem \ref{main_a} implies that, for any geometric limit $G$ of an algebraically convergent sequence $\{\rho_n\}$ of quasi-Fuchsian representations, $\mathbf{H}^3/G$ is homeomorphic to an open subset of $\Sg\times (0,1)$.
This fact is also announced by Brock, Canary and Minsky \cite{bcm}.
\end{remark}

\begin{mtheorem}\label{main_b}
For any crevasse $\cx$, there exists a geometric limit $G$ of an algebraically convergent sequence of quasi-Fuchsian representations such that $\mathbf{H}^3\cup \O(G)/G$ is homeomorphic to $\Sg\times I \setminus\cx$.
\end{mtheorem}

\begin{remark}
Though the statement of Theorem \ref{main_b} concerns only the topological type of $\mathbf{H}^3/G$, the proof given later says about some of geometric structures on $\mathbf{H}^3/G$.
In fact, one can construct $G$ so that the end $\mathcal{E}(F)$ of $\mathbf{H}^3/G$ corresponding to any front $F$ of $\cx$ has  
an arbitrarily given ending data.
Here the \emph{ending data} means  the conformal structure on the corresponding boundary component at 
infinity if $F$ is peripheral and 
the ending lamination of $\mathcal{E}(F)$ if $F$ is non-peripheral.
\end{remark}

The following problem is the geometric limit version of the Ending Lamination Theorem proved by 
Minsky \cite{mi,mi2}, Brock, Canary and Minsky \cite{bcm} and so on. 

\begin{problem}\label{P1}
Let $\cx$ be a crevasse, and let $G_i$ $(i=1,2)$ be geometric limits of algebraically convergent sequences of quasi-Fuchsian groups with homeomorphisms $h_i:\mathbf{H}^3/G_i\to \Sg\times I\setminus\cx$.
Is $h_2^{-1}\circ h_1:\mathbf{H}^3/G_1\to \mathbf{H}^3/G_2$ properly homotopic to an isometry if, for any front $F$ of $\cx$, the 
corresponding ends $\mathcal{E}_i(F)$ of $\mathbf{H}^3/G_i$ have the same ending data\,?
\end{problem}

The paper is organized as follows.
Section \ref{S_1} gives the fundamental notation and definitions needed in later sections.
In Section \ref{S_2}, topological and geometrical properties of block complexes and crevasses 
are studied.
In particular, Proposition \ref{key} shows that any infinite block admits an embedding into $\Sg\times I$.
Moreover, we prove that the complement of any crevasse in $\Sg\times I$ admits an (in)finite block decomposition.
Section \ref{S_3} shows that the geometric limit manifold $M=\mathbf{H}^3/G$ has a block complex structure, which is a key fact for the proof of Theorem \ref{main_a}.
Section \ref{S_4} presents the proof of Theorem \ref{main_a}.
Finally, in Section \ref{S_5}, we will prove Theorem \ref{main_b}.

\subsection*{Acknowledgement}
The author would like to thank Ken'ichi Ohshika for helpful comments and suggestions.

\section{Preliminaries}\label{S_1}

First of all, we will review briefly the fundamental notation and definitions needed in the paper, and refer to Thurston \cite{tl}, Benedetti and Petronio \cite{bp}, Matsuzaki and Taniguchi \cite{mt}, Canary, Epstein and Green \cite{ceg} and so on for details on hyperbolic geometry, and to Fathi, Laudenbach and Po\'{e}naru \cite{flp} for measured 
foliations (laminations).

Let $A$ be a closed subset of a metric space $(X,d)$.
For any $R>0$, the $R$-neighborhood $\{x\in X;\, d(x,A)\leq R\}$ is denoted by $\cn_R(A,X)$ (or $\cn_R(A)$ for simplicity).
In the case when $A$ is a one point set $\{x\}$, the $R$-neighborhood is also denoted by $\cb_R(x,X)$ and said to be the \emph{ball} 
in $X$ of radius $R$ centered at $x$.
The closure of $A$ in $X$ is denoted by $\ol A$.
The \emph{frontier} $\fr(A)$ of $A$ in $X$ is the intersection $\ol A\cap \ol{X-A}$.
For example, if $X=\mathbf{R}\times [0,1]$ and $A=(0,1]\times [0,1]$, then $\fr(A)=\{0,1\}\times [0,1]$.

Throughout the remainder of this paper, we suppose that $I$ is the closed interval $[0,1]$ and 
$\Sg$ is a closed orientable surface with a fixed hyperbolic structure.
A compact subsurface $F$ of $\Sg$ is \emph{geodesic} if each component of $\partial F$ is a geodesic loop in $\Sg$.
Geodesic subsurfaces (resp.\ simple geodesic loops) of $\Sg$ are simply called \emph{g-subsurfaces} (resp.\ \emph{g-loops}) of $\Sg$.
The complement in $\Sg$ of a disjoint union of g-subsurfaces and g-loops in $\Sg$ is called an \emph{open g-subsurface} of $\Sg$.
Note that, for an open g-subsurface $F$, the frontier $\fr (F)$ of $F$ in $\Sg$ may have a g-loop 
component contained in $\Int (\ol F)$ for the closure $\ol F$ of $F$ in $\Sg$.

Hyperbolic 3-space $\mathbf{H}^3$ is the Riemannian 3-manifold with the underlying space $\mathbf{C}\times \mathbf{R}_+=\{ (z,t)\in \mathbf{C}\times \mathbf{R};t>0\}$ and the metric $ds^2=(|dz|^2+dt^2)/t^2$.
One can regard that the Riemann sphere $\widehat{\mathbf{C}}=\mathbf{C}\cup \{\infty\}$ is the boundary $\part \mathbf{H}^3$ of $\mathbf{H}^3$.
The group $\mathrm{Isom}^+(\mathbf{H}^3)$ of orientation-preserving, isometric transformations on $\mathbf{H}^3$ is naturally identified with the group $\mathrm{PSL}_2(\mathbf{C})$ of M\"{o}bius transformations on $\widehat{\mathbf{C}}$.

A \emph{Kleinian group} is a discrete subgroup of $\mathrm{PSL}_2(\mathbf{C})$.
In this paper, we always assume that all Kleinian groups $\Gamma$ are torsion-free, or equivalently $\Gamma$ contains no elliptic elements.
For a Kleinian group $\Gamma$, the quotient space $M=\mathbf{H}^3/\Gamma$ is called a \emph{hyperbolic $3$-manifold}.
The \emph{limit set} $\L(\G)$ of $\G$ is the set of accumulation points of the orbit space $\G x_0$ in the 3-ball $\mathbf{H}^3\cup \part \mathbf{H}^3$ for a fixed point $x_0\in \mathbf{H}^3$.
Note that $\L(\G)$ is contained in $\part \mathbf{H}^3$.
The complement $\O(\G)=\part \mathbf{H}^3\setminus \L(\G)$ is called the \emph{domain of discontinuity} of $\G$.
One can suppose that $M$ is the interior of the manifold $\mathbf{H}^3\cup \O(\G)/\G$, which is called 
the \emph{Kleinian manifold} with the interior $M$ and denoted by $M^\natural$.
The convex hull $H_\G$ of $\L(\G)$ in $\mathbf{H}^3$ is $\G$-invariant.
The quotient $C_\G=H_\G/\G$ is called the \emph{convex core} of $M$.
By Marden \cite{ma}, $C_\G$ is homeomorphic to $M^\natural$ if $\G$ is not conjugate in $\mathrm{PSL}_2(\cc)$ to a subgroup of $\mathrm{PSL}_2(\rr)$.
The Kleinian group is \emph{geometrically finite} if the volume of the $\delta$-neighborhood of $H_\G$ in $M$ is finite.
For an $\varepsilon >0$, the $\varepsilon$-\textit{thin part} 
$M_{\mathrm{thin}(\varepsilon)}$ of $M$ is the set consisting of all points 
$x\in M$ such that there exists a non-contractible loop $l$ in $M$ with 
$l\ni x$ and of length $\leq \varepsilon$.
The complement $M_{\mathrm{thick}(\varepsilon)}=M\setminus \mathrm{Int} M_{\mathrm{thin}(\varepsilon)}$ is called the $\varepsilon$-\emph{thick part} of $M$.
A \emph{Margulis tube} is an embedded, equidistant, tubular neighborhood of a short closed geodesic in $M$.
A $\mathbf{Z}$ or a $\mathbf{Z}\times \mathbf{Z}$-\emph{cusp} $P$ is a subset of $M$ such that each component of $p^{-1}(P)$ is a horoball the stabilizer of which in $\Gamma$ is isomorphic to either $\mathbf{Z}$ or $\mathbf{Z}\times\mathbf{Z}$, where $p:\mathbf{H}^3\to M$ is the universal covering.
By the Margulis Lemma \cite[Corollary 5.10.2]{tl}, there exists a constant $\delta_0>0$ independent 
of $\G$, called a \emph{Margulis constant}, such that, for any 
$\varepsilon >0$ less than $\delta_0$, each component of 
$M_{\mathrm{thick}(\varepsilon)}$ is either a Margulis tube or $\mathbf{Z}$ or $\mathbf{Z}\times \mathbf{Z}$-cusp.
Let $M_{\mathrm{p:thin}(\varepsilon)}$ be the union of all parabolic cusp components of $M_{\mathrm{thin}(\varepsilon)}$.
We denote $M \setminus \Int M_{\mathrm{p:thin}(\varepsilon)}$ by $M_{\mathrm{p:thick}(\varepsilon)}$.

A sequence $\{\Gamma_n\}$ of Kleinian groups converges \emph{geometrically} to a subgroup $G$ of $\mathrm{PSL}_2(\mathbf{C})$ if (i) each $\gamma\in G$ is the limit of a sequence $\{\gamma_n\}$ with $\gamma_n\in \Gamma_n$ and (ii) the limit of any convergent sequence $\{\gamma_{n_i}\}$ with $\gamma_{n_i}\in \Gamma_{n_i}$ is an element of $G$.
We set $N_n=\mathbf{H}^3/\Gamma_n$ and $M_\infty=\mathbf{H}^3/G$.
A homeomorphism $g:(X,d_X)\longrightarrow (Y,d_Y)$ between metric spaces is called a $K$-\emph{bi-Lipschitz map} for a $K\geq 1$ if $g$ satisfies $d_X(x,y)/K\leq d_Y(g(x),g(y))\leq Kd_X(x,y)$ for all $x,y\in X$.
It is well known that $\{ \Gamma_n\}$ converging geometrically to $G$ is equivalent to the existence of $K_n$-bi-Lipschitz homeomorphisms $g_n:\mathcal{B}_{R_n}(x_n,N_n)\longrightarrow \mathcal{B}_{R_n}(z_\infty, M_\infty)$ with $K_n\searrow 1$, $R_n\nearrow \infty$ for a suitable choice of their base points $x_n\in N_n$, $z_\infty\in M_\infty$.
Refer to \cite{jm}, \cite[Chapter E]{bp} for more details on properties of geometric limit groups.

The holonomy $\rho:\Pi=\pi_1(M)\longrightarrow \mathrm{PSL}_2(\mathbf{C})$ of a hyperbolic 3-manifold $M$ is a faithful discrete representation.
We call that such a $\rho$ is a \emph{Kleinian representation}.
A sequence $\{\rho_n\}$ of Kleinian representations of $\Pi$ converges \emph{algebraically} to a representation $\rho_\infty:\Pi\longrightarrow \mathrm{PSL}_2(\mathbf{C})$ if, for each $\gamma\in \Pi$, $\{\rho_n(\gamma)\}$ converges to $\rho_\infty(\gamma)$ in $\mathrm{PSL}_2(\mathbf{C})$.
If $\Pi$ is \emph{non-elementary}, that is, $\Pi$ contains no abelian subgroups of finite index, then the algebraic limit $\rho_\infty$ is also a Kleinian representation (see \cite[Theorem 7.12]{mt}).
Moreover, $\{\rho_n(\Pi)\}$ has a subsequence $\{\rho_j(\Pi)\}$ converges geometrically to a Kleinian group $G\supset \rho_\infty(\Pi)$ (see \cite[Corollary 9.1.8]{tl}).
We refer to Kerckhoff and Thurston \cite{kt} and Brock \cite{br} for typical examples of geometric limits of such sequences which are not equal to their algebraic limits.
The subsequence $\{\rho_j\}$ is said to converge \emph{strongly} to $\rho_\infty$ if $\rho_\infty(\Pi)=G$.

We are above all interested in the case of $\Pi=\fd(\Sg)$.
A Kleinian representation $\rho:\Pi=\pi_1(\Sigma)\longrightarrow \mathrm{PSL}_2(\mathbf{C})$ is called \emph{quasi-Fuchsian} if the limit set of $\G=\rho(\Pi)$ is a Jordan curve in $\hat \cc$.
It is well known that, if $\rho$ is quasi-Fuchsian, then $\G$ is geometrically finite and $\mathbf{H}^3\cup \O(\G)/\G$ is homeomorphic to $\Sigma\times I$.

Throughout the remainder of this section, suppose that $F$ is either $\Sg$ itself or an open g-subsurface of $\Sg$.

A \emph{(geodesic) lamination} $\lambda$ on $F$ is a compact subset of $F$ consisting of mutually disjoint, simple geodesics, which are called \emph{leaves} of $\lambda$.
A lamination $\lambda$ is \emph{full} if $\lambda$ is not a proper sub-lamination of any other lamination on $F$.
If a lamination $\lambda$ on $F$ has an invariant transverse measure, then it is called a \emph{measured lamination}.
Let $\mathcal{ML}_0(F)$ be the set of all measured laminations on $F$, and let $\mathcal{PL}_0(F)=\bigl(\mathcal{ML}_0(F)\setminus\{0\}\bigr)/\sim$ be its projective space, where $\lambda, \lambda'\in \mathcal{ML}_0(F)\setminus\{0\}$ being equivalent, $\lambda\sim\lambda'$, means $\lambda=r\lambda'$ for some $r>0$ and the equivalence class is denoted by $[\lambda]\in \mathcal{PL}_0(F)$.
The standard topology on $\mathcal{ML}_0(F)$ and the quotient topology on $\mathcal{PL}_0(F)$ are given in \cite[Section 8.10]{tl}.
Then $\mathcal{PL}_0(F)$ is homeomorphic to the $(6h-7+2b)$-dimensional sphere $S^{6h-7+2b}$ (see \cite{tb}, \cite[Expos\'{e} 6]{flp}), where $h$ is the genus of $F$ and $b$ is the number of components of $\fr(F)$.
In particular, if $F$ is not an open pair of pants, then $\dim \mathcal{PL}_0(F)\neq 0$.
Let $\mathrm{Teich} (F)$ be the Teichm\"{u}ller space of complete hyperbolic structures on $F$ of finite area, and let $\overline{\mathrm{Teich}(F)}=\mathrm{Teich}(F)\cup \mathcal{PL}_0(F)$ be the Thurston compactification of $\mathrm{Teich}(F)$ given in \cite{tb}, see also \cite[Expos\'{e} 8]{flp}.

A $\fd$-injective proper map $f$ from $F(\sigma)$ with $\sigma\in \mathrm{Teich} (F)$ to a hyperbolic 3-manifold $M$ is called a \emph{pleated map} realizing a geodesic lamination $\lambda$ in $F(\sigma)$ if $f$ satisfies the following conditions.
\begin{enumerate}[\rm (i)]
\item $f$ maps each parabolic cusp of $F(\sigma)$ into a parabolic cusp in $M$.
\item $f$ is \emph{arcwise isometric}, that is, for any rectifiable path $\alpha$ in $F(\sigma)$, $f(\alpha)$ is also a rectifiable path in $M$ with $\mathrm{length}_{F(\sigma)}(\alpha)=\mathrm{length}_M(f(\alpha))$.
\item $f(l)$ is a geodesic in $M$ for each leaf $l$ of $\lambda$.
\item $f|\Delta$ is a totally geodesic immersion into $M$ for each component $\Delta$ of $\Sigma\setminus \lambda$.
\end{enumerate}
Then the image $f(F)$ is called a \emph{pleated surface} realizing $\lambda$.
If $T$ is an open pair of pants, then any pleated map $f:T\to M$ is a totally geodesic immersion.

Let $f:F\to M$ be a $\fd$-injective proper embedding excising from $M$ a submanifold $E$ facing an end $\mathcal{E}$ of $M$.
Consider the case when $F$ does not contain any g-loop $l$ such that $f(l)$ is 
homotopic in $E$ to a $\mathbf{Z}$-cusp of $M$.
Then we say that $\ce$ is a \emph{geometrically finite end} if there exists a collar neighborhood $N(f (F))$ of $f(F)$ in $E$ 
homeomorphic to $F\times I$ and such that $\partial N(f(F))\setminus f(F)$ is locally convex in 
$N(f(F))$.
If $\ce$ is not geometrically finite, then it is called a \emph{geometrically infinite end}.
We refer to Bonahon \cite{bo} for various results and properties on such ends.
Note that, if $\G$ is geometrically finite, then each end of $M$ is geometrically finite.

\section{Block complexes}\label{S_2}

Recall that $\mathrm{pr}_\mathrm{hz}:\Sg\times I\to \Sg$ and $\mathrm{pr}_{\mathrm{vt}}:\Sg\times I\to I$ are the projections onto the first and second factors respectively.
The preimage $\mathrm{pr}_{\mathrm{vt}}^{-1}(y)=\Sg\times \{y\}$ $(y\in I)$ is denoted by $\Sg_y$.
For any subset $R$ of either $\Sg_y$ or $\Sg$ and a subinterval $J$ (resp.\ a point $a$) of $I$, the subset $\mathrm{pr}_\mathrm{hz}(R)\times J$ (resp.\ $\mathrm{pr}_\mathrm{hz}(R)\times \{a\}$) of $\Sg\times I$ is denoted by $R_J$ (resp.\ $R_{\{a\}}$).
In particular, $R=R_{\{y\}}$ if $R\subset \Sg_y$.

First of all, we will present the definition of block complexes.
Let $F$ be either an open g-subsurface in $\Sg$ or $\Sg$ itself.
For any $0\leq a<b\leq 1$, $F_{(a,b)}$, $F_{(a,1]}$, $F_{[0,b)}$ and $F_{[0,1]}$ are called \emph{blocks} in $\Sg\times I$.
The \emph{frontier} $\fr(B)$ of a block $B=F_J$ is $F_{\part \ol{J}}$.
An open subsurfaces $F_{\{c\}}$ in $\fr(B)$ is a \emph{non-peripheral front} (resp.\ \emph{peripheral front}) of $B$ if $c\not\in J$ (resp.\ 
$c\in J$). 
From the definition, any peripheral front of $B$ is contained in $(\Sg_0\cup \Sg_1)\cap B$.
Note that, though $F_{\{0\}}$ is a subset of $\Sg_0$, $F_{\{0\}}$ is a non-peripheral front of $F_{(0,b)}$.
The union of $B\cup \fr(B)$ is denoted by $B^\bullet$.
For a block $B=F_J$, let $\alpha(B),\beta(B)$ be the elements of $\part\ol{J}$ with $\alpha(B)<\beta(B)$.
For example $\alpha(F_{[0,b)})=0$ and $\beta(F_{[0,b)})=b$.

We say that a set $\ck=\{B_0,B_1,\dots,B_n\}$ of finitely many mutually disjoint blocks is a 
\emph{block complex} if it satisfies the following conditions.
\begin{enumerate}[\rm (i)]
\item
$B_0$ contains $\Sg_{1/2}$.
\item
If $B_i^\bullet \cap B_j^\bullet\neq \eset$ for $i\neq j$, then the intersection $J_{ij}$  is a (possibly disconnected) open g-subsurface, which is called the \emph{joint} of $B_i$ and $B_j$.
\item
The union $|\ck|$ of all blocks in $\ck$ and their joints, called the \emph{support} of $\ck$, 
is connected.
\end{enumerate}

In particular, the support $|\ck|$ of any block complex is a 3-manifold containing $\Sg_{1/2}$ 
and with $\part |\ck|=|\ck|\cap (\Sg_0\cup \Sg_1)$.
Two block complexes $\ck$ and $\ck'$ are \emph{equivalent} to each other if there exists a homeomorphism $h:|\ck|\to |\ck'|$ which maps each block of $\ck$ onto a block of $\ck'$ and $h|\Sg_{1/2}$ is the identity of $\Sg_{1/2}$.
We should remark that the complements $\Sg\times I\setminus |\ck|$ and $\Sg\times I\setminus |\ck'|$ do not 
necessarily have the same topological type even if $\ck$ is equivalent to $\ck'$.

\subsection{Embeddings of infinite block complexes}
Let $\{\ck_n\}_{n=1}^\infty$ be a sequence of block complexes with $\ck_n\subsetneqq \ck_{n+1}$.
The union $\ck_\infty=\bigcup_{n=1}^\infty \ck_n$ is called an \emph{infinite block complexes} with 
support $|\ck_\infty|=\bigcup_{n=1}^\infty |\ck_n|$.

The following proposition is crucial in the proof of Theorem \ref{main_a}.

\begin{prop}\label{key}
Let $L=\bigcup_{n=1}^\infty L_n$ be a topological space with $L_n\subsetneqq L_{n+1}$ $(n=1,2,\dots)$ 
and admitting embeddings $\eta_n:L_n\to \Sg\times I$ satisfying the following conditions.
\begin{enumerate}[\rm (i)]
\item
For any $n$, the image $\eta_n(L_n)$ is the support of a block complex $\ck_n$.
\item
Each $\eta_{n+1}\circ \eta_n^{-1}:|\ck_n|\to |\ck_{n+1}|$ defines an equivalence from $\ck_n$ to a 
block subcomplex of $\ck_{n+1}$.
\end{enumerate}
Then there exists an embedding $h_\infty:L\to \Sg\times I$ such that $h_\infty(L_n)$ is the 
support of a block complex $\ck_n'$ which is equivalent to $\ck_n$ via $\eta_n\circ h_\infty^{-1}$.
In particular, $h_\infty(L)$ is the support of the infinite block complex $\ck_\infty=\bigcup_{n=1}^\infty \ck_n'$.
\end{prop}

For any integers $n,m$ with $1\leq n\leq m$, let $\ck_{m;n}$ be the subcomplex of $\ck_m$ with $|\ck_{m;n}|=\eta_m(L_n)$.
The sequence $\{\ck_{n_k;k}\}$ is called a \emph{subcomplex-subsequence} of $\{\ck_n\}$ if $\{\ck_{n_k}\}$ is a 
subsequence of $\{\ck_n\}$.
We often denote the new sequence again by $\{\ck_k\}$ for simplicity, that is, 
$\ck_k^{\mathrm{new}}=\ck_{n_k;k}^{\mathrm{old}}$ and $\eta_k^{\mathrm{new}}=\eta_{n_k}^{\mathrm{old}}|L_k$.

For any block $B_i^{(n)}$ of $\ck_n$ and any $m\geq n$, the block $\eta_m\circ \eta_n^{-1}(B_i^{(n)})$ in $\ck_{m}$ is denoted by $B_i^{(m)}$.
Set $\alpha_{i,n}=\alpha(B_i^{(n)})$ and $\beta_{i,n}=\beta(B_i^{(n)})$ for short.
Let $T_n$ be the finite subset of $I$ consisting of the numbers $\alpha_{i,n},\beta_{i,n}$ for all $B_i^{(n)}\in \ck_n$.
Consider the correspondence 
$\tau_{n,m}:T_n\to T_m$ which transfers $\alpha_{i,n},\beta_{i,n}$ respectively to $\alpha_{i,m},\beta_{i,m}$.
Note that $\tau_{n,m}$ is not necessarily a map.
In fact, it may occur that  $\alpha_{i,n}=\alpha_{j,n}$ (resp.\ $\alpha_{i,n}=\beta_{j,n}$) but $\alpha_{i,m}\neq \alpha_{j,m}$ 
(resp.\ $\alpha_{i,m}\neq \beta_{j,m}$) and so on.

For the proof, we need the following two rearrangements on $\ck_n$.

\subsection*{Rearrangement I}
If necessary passing to a subcomplex-subsequence of $\{\ck_n\}$, we may assume that for any $n,m$ with $n\leq m$, the correspondence 
$\tau_{n,m}:T_n\to T_m$ is an order-preserving map.
Deform the new $\ck_n$'s by horizontal ambient isotopies of $\Sg\times I$, we may assume that 
$\alpha_{i,n}=\alpha_{i,m},\beta_{i,n}=\beta_{i,m}$ for any $n,m$ with $n\leq m$ and $i$ with 
$B_i^{(n)}\in \ck_n$.
In particular, $T_n\subset T_m$.

\subsection*{Rearrangement II}
\begin{figure}[hbtp]
\centering
\scalebox{0.5}{\includegraphics[clip]{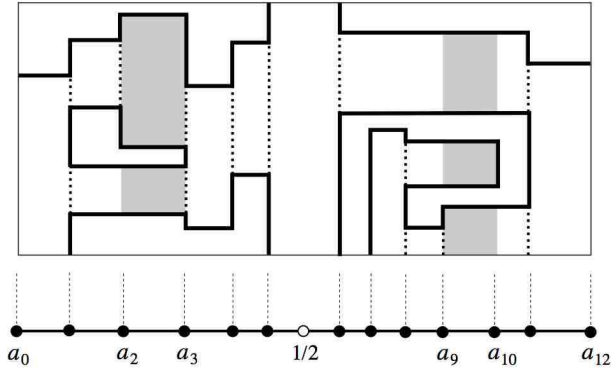}}
\caption{The union of the shaded regions in the left (resp.\ right) hand side 
is $\eta_n(R_3^n)$ (resp.\ $\eta_n(R_{10}^n)$).}
\label{rear_2}
\end{figure}

Set $T_n=\{a_0,a_1,\dots,a_t\}$ with $a_{j-1}<a_j$,  $P_j^n=\Sg_{[a_{j-1},a_j]}$ for $j=1,\dots,t$ and $R_j^n=\eta_n^{-1}(P_j^n)$, see Fig.\ \ref{rear_2}.
If necessary passing again to a subcomplex-subsequence of $\{\ck_n\}$, we may assume that all $\eta_m|R_j^n$ $(m\geq n)$ define 
\emph{the same embedding up to marking}, that is, there exists an orientation-preserving 
homeomorphism $\gamma_{m,n}:P_j^n\to P_j^n$ with 
$\gamma_{m,n}\circ \eta_m|R_j^n=\eta_{n}|R_j^n$ for all $m\geq n$.

\bigskip

\begin{figure}[hbtp]
\centering
\scalebox{0.5}{\includegraphics[clip]{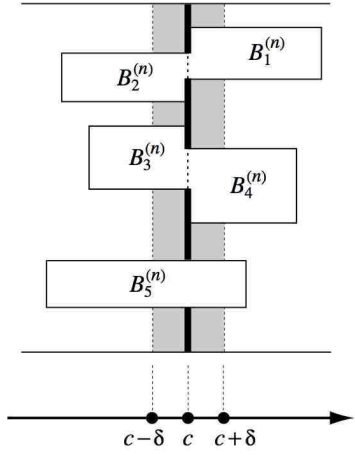}}
\caption{The union of the bold vertical segments represents $S_c^{(n)}$.
The union of $S_c^{(n)}$ and the shaded regions is $Q_\delta(S_c^{(n)})$.}
\label{d_region}
\end{figure}

We say that $S_c^{(n)}=\Sg_c\setminus |\ck_n|$ with $c\in I$ is the \emph{slit} for $|\ck_n|$ at $c$.
By Rearrangements I and II, for all sufficiently large $n$, $S_c^{(n)}$ have the same topological type.
Set $-\chi_{\min}(S_c^{(n)})=\min_{m\geq n}\{-\chi(S_c^{(m)})\}$.
The slit $S_c^{(n)}$ is \emph{stable} if all $S_c^{(m)}$ $(m\geq n)$ are homeomorphic to each other or 
equivalently $-\chi(S_c^{(n)})=-\chi_{\min}(S_c^{(n)})$.
Let $T_\infty'$ be the set of accumulation points of $T_\infty=\bigcup_{n\geq 1}T_n$.
For $c\in T_\infty'$, suppose that $B_1^{(n)},\dots,B_k^{(n)}$ are the blocks in $\ck_n$ with $B_i^{(n)\bullet}\cap \Sg_c\neq \eset$.
Take $\delta>0$ sufficiently small so that $\Sg_{[c-\delta,c)}\cup \Sg_{(c,c+\delta]}$ does not meet any fronts 
of $B_i^{(n)}$ $(i=1,\dots,k)$.
Then the set 
$$Q_\delta(S_c^{(n)})=\bigl(\Sg_{[c-\delta,c)}\cup \Sg_{(c,c+\delta]}\setminus B_1^{(n)}\cup\dots\cup B_k^{(n)}\bigr)\cup S_c^{(n)}$$
is called the \emph{$\delta$-region} of the slit $S_c^{(n)}$ for $|\ck_n|$, see Fig.\ \ref{d_region}.
For any integer $s\geq 1$, let $T_{\infty,s}'$ be the subset of $T_\infty'$ consisting of elements $c\in T_\infty'$ with 
$-\chi_{\min}(S_c^{(n)})=s$.
Suppose that $S_c^{(n)}$ is stable.
If $S_d^{(n)}$ $(d\in T_\infty')$ is contained in $Q_\delta(S_c^{(n)})\setminus S_c^{(n)}$ and satisfies 
$-\chi_{\min}(S_d^{(n)})=s$, then $S_d^{(n)}$ is parallel to $S_c^{(n)}$ in $\Sg\times I\setminus |\ck_n|$ for all sufficiently large $n$.
Thus there exists $\delta(c)>0$ for $c\in T_{\infty,s}'$ independent of $n$ such that $-\chi_{\min}(S_d^{(n)})<
s$ if $S_d^{(n)}\subset 
Q_{\delta(c)}(S_c^{(n)})\setminus S_c^{(n)}$.
In particular, this implies that all $T_{\infty,s}'$ and hence $T_\infty'$ are countable subsets of $I$.

\begin{figure}[hbtp]
\centering
\scalebox{0.5}{\includegraphics[clip]{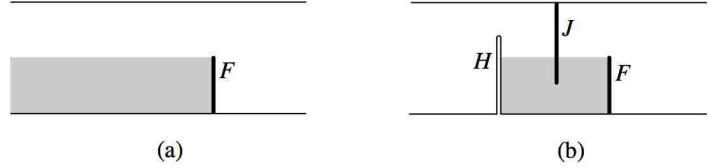}}
\caption{The shaded parts represent the affected regions.}
\label{sliding}
\end{figure}

Let $F$ be a g-subsurface of $\Sg_a$ with $0<a<1$ and $\varphi:F\to F$ an orientation-preserving homeomorphism 
such that $\varphi|\part F$ is the identity of $\part F$.
Consider the 3-manifold $M_\varphi$ obtained from $\Sg\times I\setminus \Int F$ by identifying the $(\pm)$-sides 
$F^{(\pm)}$ of $F$ via $\varphi:F^{(-)}\to F^{(+)}$.
Then $\Sg\times I\setminus \Int F$ is naturally regarded as a subset of $M_\varphi$. 
We say that $M_\varphi$ is the manifold obtained from $\Sg\times I\setminus \Int F$ by the \emph{$\varphi$-sliding along} $F$.
There exists a homeomorphism $\xi_0:M_\varphi\to \Sg\times I$ such that the restriction $\xi_0|(\Sg\times I\setminus 
C_0)$ is the identify, where $C_0$ is either $F_{[0,a)}$ or $F_{(a,1]}$.
The $C_0$ is called the \emph{affected region} of the sliding, see Fig.\ \ref{sliding}\,(a).
For the proof of Proposition \ref{key}, we reduce the affected region by 
adopting the following minor trick.
Let $H$ be a non-peripheral vertical g-subsurface in $\Sg\times I$ not in $\Sg_a$, say $H\subset \Sg_b$ for some 
$0<b<a$, with $\mathrm{pr}_\mathrm{hz}(H)\supset \mathrm{pr}_\mathrm{hz}(F)$.
Then there exists a homeomorphism $\xi_1: M_\varphi\setminus H\to \Sg\times I\setminus H$ such the the affected 
region is $C_1=F_{(b,a)}$, see Fig.\ \ref{sliding}\,(b).
In the case when there exists another g-subsurface $J$ of $\Sg_c$ for some $b<c<a$ satisfying $\mathrm{pr}_\mathrm{hz}(\Int F)
\cap \mathrm{pr}_\mathrm{hz}(\Int J)\neq \eset$ and $\mathrm{pr}_\mathrm{hz}(F)\setminus \mathrm{pr}_\mathrm{hz}(J)\neq \eset$.
Then we do not count on any partial contribution of $J$ to the reduction of the affected region of the $\varphi$-sliding.
More precisely, the homeomorphism $\xi_2:M_\varphi\setminus H\cup J\to \Sg\times I\setminus H\cup J'$ is the restriction 
of $\xi_1$, where $J'=\xi_1(J)$.
In particular, the affected region of $\xi_2$ is $C_1\setminus J$. 
If the affected region is contained a subset $Q$ of $\Sg\times I$, then we say that 
the effect of the sliding is \emph{absorbed} in $Q$.

\begin{proof}[Proof of Proposition \ref{key}]
For any $n\geq 2$, $\eta_{n-1}:L_{n-1}\to \Sg\times I$ is an embedding with 
$\part L_{n-1}=\eta_{n-1}^{-1}(\Sg_0\cup \Sg_1)$, but in general $\eta_n|L_{n-1}$ is not 
equal to $\eta_{n-1}$.
We set $h_{n-1}=\eta_{n-1}$ for some $n\geq 2$ and $T_n=\{a_0,a_1,\dots, a_t\}$ as above.
By Rearrangement II, for any $j=1,\dots,t$, there exists an embedding $\hat h_n^j:R_n^j\to P_n^j$ such that $\hat h_n^j|R_n^j\cap L_{n-1} 
= h_{n-1}|R_n^j\cap L_{n-1}$.
Note that the union of $\hat h_n^j$ does not necessarily define an entire embedding from $L_n$ to $\Sg\times I$.
Let $\hat T_n$ be the subset of $T_n$ consisting of elements $a_j\in T_n$ with $-\chi(S_{a_j}^{(n-1)})>-\chi(S_{a_j}^{(n)})$, 
where $S_{a_j}^{(n-1)}=\Sg_{a_j}\setminus |\ck_{n-1}|$ and $S_{a_j}^{(n)}=\Sg_{a_j}\setminus |\ck_{n}|$.
In particular, $c\in \hat T_n$ implies that $S_c^{(n-1)}$ is unstable.
There exist orientation-preserving homeomorphisms $\varphi_{a_j}:S_{a_j}^{(n-1)}\to S_{a_j}^{(n-1)}$ $(a_j\in \hat T_n)$ such that $\varphi_{a_j}|\part S_{a_j}^{(n-1)}$ is the identity of $\part S_{a_j}^{(n-1)}$ and 
$\bigcup_{j=1}^t \hat h^j_n$ extends to an embedding $\hat h_n:L_n\to M_n$ with 
$\hat h_n|L_{n-1}=h_{n-1}$, where $M_n$ is the manifold obtained from $\Sg\times I\setminus \bigcup_{a_j\in \hat T_n}S_{a_j}^{(n-1)}$ by the $\varphi_{a_j}$-slidings.
Let $\mathcal{S}_\infty^{(n)}$ be the union of stable slits $S_c^{(n)}$ with $c\in T_\infty'\setminus \hat T_n$.
Then $\mathcal{S}_\infty^{(n)}\subset \Sg\times I\setminus \bigcup_{a_j\in \hat T_n}S_{a_j}^{(n-1)}$ is naturally regarded as a subset of $M_n$.
We have a homeomorphism $\xi_n:M_n\setminus \mathcal{S}_\infty^{(n)}\to \Sg\times I\setminus {\mathcal{S}_\infty^{(n)}}'$ which is the identity on the complements of the union of the affected regions of these $\varphi_{a_j}$-slidings, where ${\mathcal{S}_\infty^{(n)}}'$ is the union of vertical g-subsurfaces in $\Sg\times I$ corresponding to the slits in $\mathcal{S}_\infty^{(n)}$.
Then $h_n=\xi_n\circ \hat h_n:L_n\to \Sg\times I$ is an embedding such that $h_n(L_n)$ is the support of a block complex 
equivalent to $\ck_n$, but in general $h_n$ is no longer an extension of $h_{n-1}$.
For all $m\geq n$, an embedding $h_m$ is defined similarly.
We will see that $\{h_n\}$ defines the same map eventually.

Take any subset $B$ of $L_m$ with $B^{(m)}=\eta_m(B)\in \ck_m$.
For any $c\in T_\infty'$, one can take $\delta(c)>0$ sufficiently small and suitably 
so that (i) $Q_{\delta(c)}(S_c^{(w)})$ is disjoint from 
$B^{(w)}=\eta_{w}(B)$ for all $w\geq m$ and (ii) either $[c-\delta(c),c+\delta(c)]$ 
and $[c'-\delta(c'),c'+\delta(c')]$ are mutually disjoint or one of them contains the other for any $c,c'
\in T_\infty'$.
Since $T_\infty\cup T_\infty'$ is compact, there exists a finite subset $\{c_1,\dots,c_k\}$ of 
$T_\infty'$ such that $T_\infty\setminus \bigcup_{i=1}^k[c_i-\delta(c_i),c_i+\delta(c_i)]$ contains 
only finitely many elements $b_1,\dots,b_u$, see Fig.\ \ref{sliding2}.
\begin{figure}[hbtp]
\centering
\scalebox{0.5}{\includegraphics[clip]{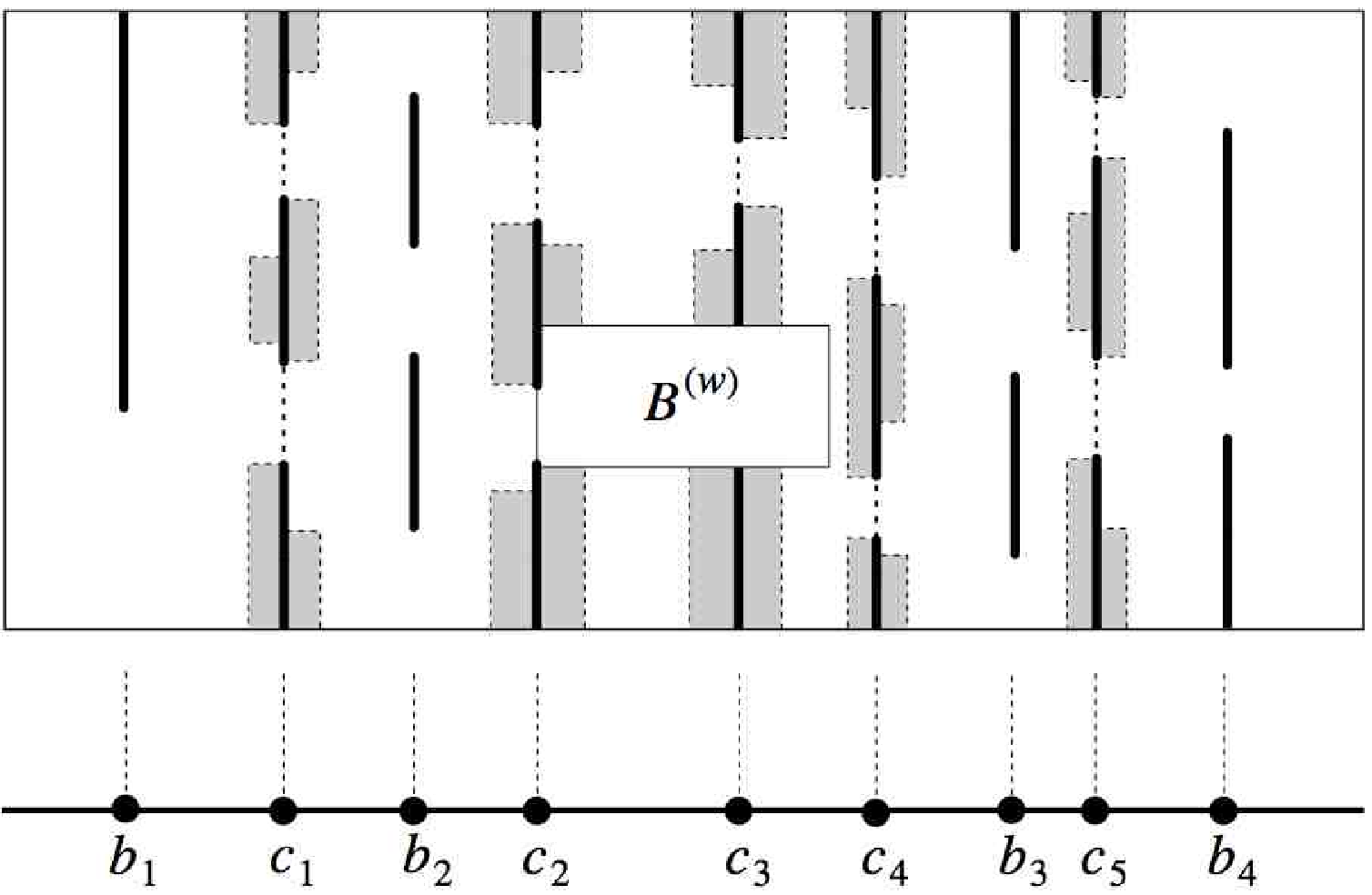}}
\caption{}
\label{sliding2}
\end{figure}
Take $w_0\in \nn$ such that $S_j^{(w_0)}$ is stable for all $j\in \{b_1,\dots,b_u,c_1,\dots,c_k\}$, that is, 
all slidings along slits at $b_1,\dots,b_u,c_1,\dots,c_k$ are already finished until the $w_0$-th step.
Since the effect of any sliding along the slit $S_d^{(w)}$ at $d\in [c_i-\delta(c_i),c_i+\delta(c_i)]$  
is absorbed in $Q_{\delta(c_i)}(S_{c_i}^{(w)})$ and $B^{(w)}\cap 
Q_{\delta(c_i)}(S_{c_i}^{(w)})=\eset$ for $i=1,\dots,k$ and $w\geq w_0$, it follows that all $h_{w}|B$ with 
$w\geq w_0$ are the same map.
Thus our desired embedding $h_\infty:L\to \Sg\times I$ is well defined by $h_\infty|B=h_{w_0}|B$.
\end{proof}

\subsection{Properties of crevasses}\label{p_c}

Here we  study topological properties of a crevasse $\cx$, that is, $\cx$ is a closed subset of $\Sg\times I$ satisfying 
the conditions (i)-(iii) in Theorem \ref{main_a}.

Set $\cy=\mathrm{pr}_{\mathrm{vt}}(\cx)$ and $X_y=\cx\cap \Sg_y$ ($y \in \cy$).
Since $\cx$ is closed, both $\L_{y}^+=\Sg_y\cap \ol{\Sg_{(y,1]}\cap \cx}$ and $\L_{y}^-=\Sg_y\cap \ol{\Sg_{[0, y)}\cap \cx}$ are contained in $X_y$.
For a g-loop $l$ of $\Sg$, $l_{[z,w]}$ with $0\leq z<w\leq 1$ is called a \emph{horizontal annulus}. 
A g-loop $m_y\subset X_y$ is \emph{parallel} to a g-loop $m_z\subset X_z$ in $\Sg\times I\setminus \cx$ if there exists a continuous map $h$ from an annulus $A$ to $\Sg\times I$ with $h(\part A)=m_y\cup m_z$ and $h(\Int A)\subset \Sg\times I\setminus\cx$.
These loops are \emph{horizontally parallel} in $\cx$ if there exists a horizontal annulus $A'$ contained in $\cx$ with $\part A'=m_y\cup m_z$.

\begin{figure}[hbtp]
\centering
\scalebox{0.5}{\includegraphics[clip]{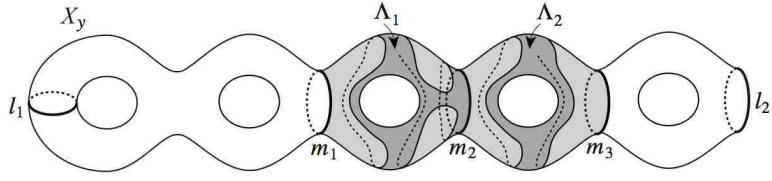}}
\caption{The case of $\L_y^\ve=l_1\cup \L_1\cup \L_2$ and $\fr(X_y)=l_2$.
The union of the two white parts is $Z(\L_y^\ve)$ and $\lambda(\L_y^\ve)=l_1\cup l_2\cup m_1
\cup m_2\cup m_3$.}
\label{d_l}
\end{figure}

For $y\in \mathcal{Y}$, let $Z(\L_y^\ve)$ be a maximal open g-subsurface of $X_y$ disjoint from $\L_y^\ve$.
A component of $Z(\L_y^\ve)$ is called a \emph{non-peripheral front} of $\cx$, and a component of 
$\Sg_0\cup \Sg_1\setminus \cx$ is a \emph{peripheral front} of $\cx$.
As was defined in Introduction, $\lambda(\L_y^\ve)$ is the union of g-loops in $X_y\setminus 
Z_y(\L_y^\ve)$ each of which has a free side with respect to $\L_y^\ve$.
The definition implies that $\lambda(\L_y^\ve)\supset \fr(X_y)$, see Fig.\ \ref{d_l}.

In Lemmas \ref{l21}-\ref{l23}, we only consider $\L^+_y$.
It will be seen easily that the corresponding assertions for $\L^-_y$ hold.

\begin{figure}[hbtp]
\centering
\scalebox{0.5}{\includegraphics[clip]{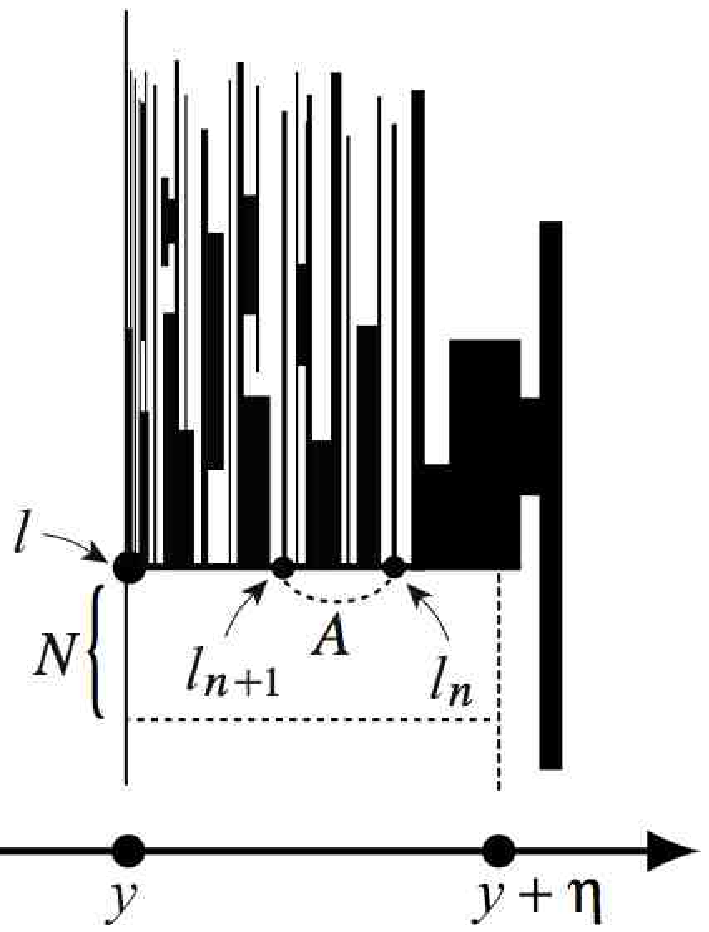}}
\caption{The black part represents $\cx$.}
\label{collar}
\end{figure}

\begin{lemma}\label{l21}
If $l$ is a g-loop in $\lambda(\L_y^+)$ with $l\cap \L_y^+\neq \emptyset$, then there exists a collar neighborhood $N$ of $l$ in $\Sg_y$ and a constant $\eta=\eta(y)>0$ such that $N_{(y,y+\eta]}\cap \cx=l_{(y,y+\eta]}$.
In particular, $l$ is contained in $\L_y^+$, e.g.\ $m_2$ in Fig.\ \ref{d_l}.
\end{lemma}
\proof
Since $l$ has a free side with respect to $\L_y^+$, for a sufficiently small $\delta >0$, at least one component $N'$ of $\cn_\delta (l,\Sg_y)\setminus l$ is disjoint from $\L_y^+$.
Then $N=N'\cup l$ is an annulus with $\part N=\partial N'\cup l$.
For any $x\in l\cap \L_y^+$, there exists a sequence $\{x_n\}$ with $x_n\in X_{y_n}$  $(y_n\searrow y)$ converging to $x$.
Let $R_n$ be the component of $X_{y_n}$ containing $x_n$.
Since $N'\cap \L_y^+=\eset$, the intersection $\partial N'\cap (R_{n})_{\{y\}}$ is empty for all sufficiently large $n\in \nn$.
Since $l$ is a deformation retract of the annulus $N$, $N\cap (R_{n})_{\{y\}}\subset l$.
Thus we have $n_0\in \nn$ such that $N_{(y,y_{n_0}]}'\cap \cx=\emptyset$.
Since $\{(x_{n})_{\{y\}}\}\sto x$ in $\Sg_y$, $\partial R_n$ has a component $l_n$ with $(l_{n})_{\{y\}}=l$.
We set $\eta=y_{n_0}-y$.
For any $n\geq n_0$, there exists an embedded annulus $A$ in $N_{(y,y+\eta]}$ with $\Int A\subset N'_{(y,y+\eta]}$ and $\partial A=A\cap l_{(y,y+\eta]}=l_n\cup l_{n+1}$, see Fig.\ \ref{collar}.
By the condition (iii), $l_{[y_n,y_{n+1}]}$ is contained in $\cx$.
From this, we know that $l_{(y,y+\eta]}\subset \cx$ and hence $N_{(y,y+\eta]}\cap \cx=l_{(y,y+\eta]}$.
\qed\\

\begin{remark}\label{R2}
In the proof of Lemma \ref{l21}, we only use the weak version of the condition (iii) such that, for any $y,z\in \cy$ with $y<z$, if a component $l_y$ of $\fr (X_y)$ is parallel to a component $l_z$ of $\fr (X_z)$ in $\Sg\times I \setminus\cx$, then $l_y$ and $l_z$ are horizontally parallel in $\cx$.
This fact is crucial in the proof of Theorem \ref{main_a}.
\end{remark}

A disjoint union of g-loops and g-subsurfaces in $\Sg_y$ is called a \emph{g-l-subsurface} 
of $\Sg_y$.
Note that $\D(\L_y^+)=X_y\setminus Z(\L_y^+)$ is a minimal g-l-subsurface of $X_y$ containing 
$\L_y^{+\bullet}=\L_y^+\cup\fr(X_y)$.

\begin{lemma}\label{l22}
$Z(\L_y^+)$ is determined uniquely.
\end{lemma}
\proof
Let $\D_1$, $\D_2$ be minimal g-l-subsurfaces of $X_y$ containing $\L_y^{+\bullet}$.
Suppose that $\D_1\neq \D_2$.
Then some component $l_1$ of $\fr(\D_1)$ meets a component $l_2$ of $\fr(\D_2)$ non-trivially and transversely, otherwise $\D_1\cap \D_2$ would be a g-l-subsurface of $\Sg_y$ containing 
$\L_y^{+\bullet}$ with $\D_1\cap \D_2\subsetneqq \D_i$ for at least one $i$ of $\{1,2\}$.
If $l_i\cap \L_y^{+\bullet}\neq \eset$ for an $i\in \{1,2\}$, then by Lemma \ref{l21} $l_i$ would be contained in $\L_y^{+\bullet}\subset \D_{3-i}$, a contradiction.
Thus $\D_1\setminus l_1\cup l_2$ contains $\L_y^{+\bullet}$.
Let $\D_1'$ be the union of all components of $\D_1\setminus l_1\cup l_2$ meeting $\L_y^{+\bullet}$ non-trivially, and $\hat \D_1$ a (uniquely determined) g-l-subsurface in $\D_1'$ such that each component of $\D_1'\setminus\hat \D_1$ is an open annulus.
As in the proof of Lemma \ref{l21}, there exists $\eta>0$ such that, for any $z\in \cy$ with $0<z-y\leq \eta$, $(X_z)_{\{y\}}$ is contained in $\D_1'$.
Since $\hat \D_1$ is a deformation retract of $\D_1'$ with $\fr(\hat \D_1)$ geodesic, we have $(X_z)_{\{y\}}\subset \hat \D_1$.
It follows that $\L_y^{+\bullet}$ is contained in $\hat \D_1$.
This contradicts the minimality of $\D_1$, and hence $\D_1=\D_2$.
\qed\\

\begin{lemma}\label{l23}
$\lambda(\L_y^+)$ is a disjoint union of g-loops in $\Sg_y$.
\end{lemma}
\proof
This lemma is proved by an argument quite similar to that in Lemma \ref{l22}.
If $l_1,l_2$ were g-loops in $\lambda(\L_y^+)$ meeting non-trivially and transversely, then $\D(\L_y^+)\setminus l_1\cup l_2$ would contain a g-l-subsurface $\hat \D$ with $\hat \D\supset \L_y^{+\bullet}$, a contradiction.
\qed\\

\subsection{Block complexes and crevasses}\label{ss_2_3}

Let $\cx$ be any crevasse in $\Sg\times I$.
Set
$$L=\Sg\times I\setminus \cx.$$
For any $a\in I$, a subsurface of $\Sg_a$ is called a \emph{vertical} surface in $\Sg\times I$.
Suppose that $\cf=\cf_\cx$ is the set of vertical open g-subsurfaces $F$ in $L$ with $\fr (F)\subset \cx$ and fibers $\Sg_y$ disjoint from $\cx$.
Then the union $\bigcup_{F\in \cf}\fr (F)$ consists of mutually disjoint 
horizontal annuli and g-loops.
Let $\cp_\cx$ be the set of all such annuli and loops.
An element $y\in \cy\setminus \{0,1\}$ is \emph{irregular} if $\Sg_y$ either contains a g-loop element of $\cp_\cx$ or meets two annulus components of $\cp_\cx$ one of which is contained in $\Sg_{[0,y]}$ and the other in $\Sg_{[y,1]}$.
The subset $\cy'$ of $\cy$ consisting of all irregular elements and the elements in $\{0,1\}\cap \cy$ is a countable set, $\cy'=\{y_1,y_2,\dots\}$.
Replacing $(\Sg_{y_n},X_{y_n})$ in $(\Sg\times I,\cx)$ by the saturated pair $(\Sg_{y_n}\times [0,2^{-n}],X_{y_n}\times [0,2^{-n}])$ for each $y_n\in \cy'$, 
we have a new crevasse $\cx'$ in $\Sg\times [0,2]$ such that $\Sg\times [0,2]\setminus \cx'$ is homeomorphic to $L$.
Squeezing $\Sg\times [0,2]$ in half, we may assume that $\cx'$ is a crevasse in $\Sg\times I$ and call it a \emph{saturated crevasse} obtained from $\cx$.
Note that (i) each element of $\cp_{\cx'}$ consists of horizontal annuli,  (ii) for any  $F\in \cf_{\cx'}$, there exists another $F'\in \cf_{\cx'}$ 
properly isotopic to $F$ in $L$ and (iii) any non-peripheral front of $\cx'$ is not in $\Sg_0\cup \Sg_1$.
An element $H\in \cf_{\cx'}$ is \emph{critical} if $H\cap (\Sg_0\cup \Sg_1)=\eset$ and at least one of components of $\fr(H)$ is an 
edge of an annulus in $\cp_{\cx'}$.

\begin{lemma}\label{c_to_b}
For any crevasse $\cx$, there exists a crevasse $\hat \cx$ satisfying that 
$\hat L=\Sg\times I\setminus \hat\cx$ is homeomorphic to $L$ and admits a block decomposition $\hat\ck$ such that any maximal open 
horizontal annulus in $\hat L$ meets only finitely many blocks in $\hat\ck$.
\end{lemma}
\begin{proof}
We may assume that $\cx$ is a saturated crevasse.
Let $\cf_{\mathrm{n.c.}}$ be the subset of $\cf=\cf_{\cx}$ consisting of non-critical elements.
For any $F\in \cf_{\mathrm{n.c.}}$, let $B_F$ be the union of elements in $\cf_{\mathrm{n.c.}}$ properly isotopic to $F$ in $L$.
Note that $B_F$ is a block in $\Sg\times I$ for any $F\in \cf_{\mathrm{n.c.}}$ and either $B_F=B_{F'}$ or $B_F\cap B_{F'}=\eset$ holds for any $F,F'\in \cf_{\mathrm{n.c.}}$.
Let $\ck=\{B_i\}$ be the set of such blocks.
Since $\cx$ is saturated, any non-empty $H=B_i^\bullet \cap B_j^\bullet$ with $i\neq j$ is contained in $L$, otherwise $y\in \cy$ with $\Sg_y\supset H$ would be an irregular element.
However a component $J$ of $L\setminus \bigcup \ck$ is not necessarily a joint of two blocks in $\ck$.
In fact, any small neighborhood of $J$ in $L$ may intersect infinitely many blocks in $\ck$.
Since $\cx$ is saturated, $J$ is approached by blocks in $\ck$ only from one side, say the $(+)$-side, see Fig.\ \ref{x_hat1}\,(a).
\begin{figure}[hbtp]
\centering
\scalebox{0.5}{\includegraphics[clip]{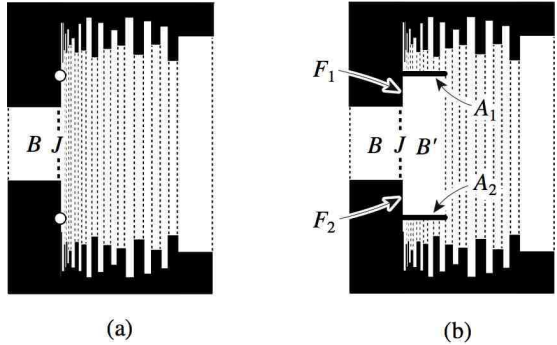}}
\caption{(a) The white dots represent $\lambda(\Lambda_y^+)$.
$B$ is a block of $\ck$ with front $J$.
(b) $A_1,A_2$ are annuli added to $\cx$ and $B'$ is a block in $\ck_J$.
$F_1,F_2$ are fronts of $\cx_J$.}
\label{x_hat1}
\end{figure}
One can construct a new crevasse $\cx_J$  by attaching narrow horizontal annuli to $\cx$ along $\lambda(\L_y^+)$ so that $J$ is a joint of two blocks in the block set $\ck_J$ associated to $\cx_J$, where 
$y\in \cy$ is the element with $\Sg_y\supset J$.
See Fig.\ \ref{x_hat1}\,(b).
\begin{figure}[hbtp]
\centering
\scalebox{0.5}{\includegraphics[clip]{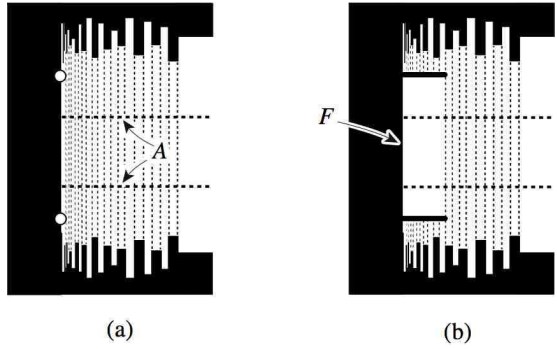}}
\caption{(a) $A$ is a part of a maximal horizontal open annulus in $L_1$ meeting infinitely many blocks in $\ck_1$.
(b) $F$ is a front of $\hat \cx$.}
\label{x_hat2}
\end{figure}
Performing similar modifications on $\cx$ for all components of $L\setminus \bigcup \ck$, we have a crevasse $\cx_1$ admitting a block complex $\ck_1$.
Modifying $\cx_1$ in similar manners in small neighborhoods of all fronts of $\cx_1$ which are  
approached by infinite sequences of blocks in $\ck_1$ (see Fig.\ \ref{x_hat2}\,(a)), we have our desired crevasse $\hat \cx$ (see Fig.\ \ref{x_hat2}\,(b)).
\end{proof}

\section{Block decomposition of geometric limits}\label{S_3}

\begin{figure}[hbtp]
\centering
\scalebox{0.5}{\includegraphics[clip]{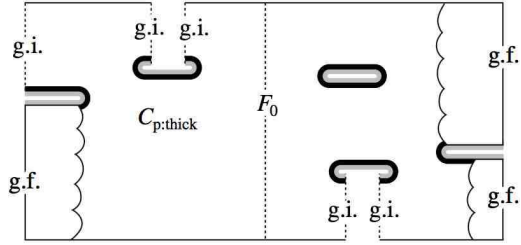}}
\caption{The union of bold curves and loops is $P$.
The union of waved curves is $\part C\setminus M_{\pthin}$.
Each segment labelled with `g.i.' (resp.\ `g.f.') represents a geometrically infinite (resp.\ finite) end of 
$M$.}
\label{c_pthick}
\end{figure}

Let $\{\rho_n:\Pi\to \mathrm{PSL}_2(\cc)\}$ be a sequence of quasi-Fuchsian representationes converging algebraically to a 
Kleinian representation $\rho_\infty$ and such that $\{\G_n\}$ converges geometrically to a Kleinian group $G$, where 
$\Pi=\fd(\Sg)$ and $\G_n=\rho_n(\Pi)$.
Note that $\G_\infty=\rho_\infty(\Pi)$ is a subgroup of $G$.
Set $N_n=\mathbf{H}^3/\G_n, N_\infty=\mathbf{H}^3/\G_\infty$ and $M=\mathbf{H}^3/G$ and denote the convex core of $N_n$ (resp.\ $M$) by $C_n$ (resp.\ $C$).
Let $p:N_\infty\to M$ be the covering associated to $\G_\infty\subset G$.
One can take an element $\mu$ of $\mathcal{ML}(\Sg)$ so that there exists a pleated map $f_n:\Sg\to C_n$ realizing $\mu$ with \emph{right marking}, that is, $\fd(f_n)=\rho_n$ for all $n\in \nn\cup \{\infty\}$.
By taking the base points $x_n$ of $N_n$ and $x_\infty\in M$ suitably, we have $K_n$-bi-Lipschitz maps $g_n:\cb_{R_n}(x_n,N_n)\to 
\cb_{R_n}(x_\infty,M)$ with $R_n\nearrow \infty$ and $K_n\searrow 1$ and such that $\{g_n\circ f_n\}$ converges uniformly to 
$p\circ f_\infty$ in $M$.
It is well know that $p\circ f_\infty$ is not necessarily homotopic to an embedding in $M$, for example see 
Anderson and Canary \cite{ac0}.
However, since $f_n$ is homotopic to an embedding into $C_n$ for any $n\in \nn$, there exists a $\fd$-injective embedding 
$h:\Sg\to C$ contained in an arbitrarily small neighborhood of $\mathrm{Image}(g_n\circ f_n)$.
We suppose that $F_0=h(\Sg)$ is the base surface of $C$ and $M$.
Take a $\delta>0$ less than a Margulis constant and such that $F_0\cap M_{\mathrm{p:thin}(\delta)}=\eset$.
From now on, for any hyperbolic 3-manifold $N$, set $N_{\mathrm{thin}(\delta)}=N_{\mathrm{thin}}$ and $N_{\mathrm{thick}(\delta)}=N_{\mathrm{thick}}$ for short.
If necessary modifying $g_n$ slightly, we may assume that $g_n(N_{n,\mathrm{thin}}\cap \cb_{R_n}(x_n,N_n))=M_{\mathrm{thin}}\cap \cb_{R_n}(x_\infty,M)$.
Though $N_n$ contains no parabolic cusps, we need to consider components of $N_{n,\mathrm{thin}}$ corresponding to 
parabolic cusps of $M$.
The union of components of $N_{n,\mathrm{thin}}$ meeting $g_n^{-1}(M_{\pthin})$ non-trivially is denoted by $N_{n,\ppthin}$.
Set $N_n\setminus \Int N_{n,\ppthin}=N_{n,\ppthick}$.
The intersections $C_n\cap N_{n,\ppthick}$ and $C\cap M_{\pthick}$ are denoted by $C_{n,\ppthick}$ and $C_{\pthick}$ 
respectively.
The convex hulls $H_{\G_n}$ of the limit sets $\L(\G_n)$ converge geometrically to the convex hull $H_G$ of $\L(G)$ (see \cite[Corollary 2.2]{kt}).
Set $P_n=P_n(\delta)=C_n\cap \part N_{n,\ppthin}$ and $P=P(\delta)=C\cap \part M_{\pthin}$ and call the pairs $(C_{n,\ppthick},P_n)$ 
and $(C_{\pthick},P)$ the \textit{pared manifolds} in $N_{n,\ppthick}$ and $M_{\pthick}$ respectively.
Furthermore, $P_n,P$ are called the \emph{parabolic loci} of these pared manifolds.
Note that the complement $\part C_{\pthick}\setminus P$ is equal to $\part C\setminus M_{\pthin}$, see Fig.\ \ref{c_pthick}.
If necessary modifying $g_n$ slightly again, we may assume that $g_n^{-1}$ maps $\part C\setminus P$ into 
$\part C_n\setminus P_n$ for all sufficiently large $n$.
In fact, $g_n^{-1}(\part C\setminus P_n)$ is the union of components of $\part C_n\setminus P_n$ 
with uniformly bounded distance from $x_n$ in $N_n$.

\subsection{Ubiquity of pseudo-pleated surfaces}

We may assume that any pleated map $f:\Sg(\sigma_n)\to N_n$ meets $\partial N_{n,\thick}$ transversely if necessary deforming $f$ slightly by a homotopy.
Let $A$ be a component of $f^{-1}(N_{n,\thin})$.
Suppose that each component of $\part A$ is non-contractible in $\Sg$.
Since $f|A:A\to N_{n,\thin}$ is $\fd$-injective and $\fd(N_{n,\thin})$ is 
abelian, $A$ is an annulus.
Since $N_{n,\thick}$ is locally convex, if a component $l$ of $\partial A$ is contractible in $\Sigma$, then $A$ is a disk bounded by $l$.
Let $\mathcal{D}_f$ and $\mathcal{A}_f$ be the union of disk and annulus components of $f^{-1}(N_{n,\thin})$ respectively.
Since the closure $\ol{\Sg\setminus \mathcal{A}_f}$ is contained in $\Sg(\sigma_n)_{\thick}$, the diameter of each component of 
$\ol{\Sg\setminus \mathcal{A}_f}$ is bounded by a constant independent of $n$.
Since $\mathrm{Area}(\mathcal{D}_f)< \mathrm{Area}(\Sigma)= -2\pi \chi(\Sigma)$, 
there exists a constant $0<\ve<\delta$ depending only on $\delta$ and $\chi(\Sg)$ such that $f(\mathcal{D}_f)\cap 
\cn_1(N_{n,\mathrm{thin}(\ve)})=\eset$.
This fact implies that $f$ can be deformed to a map $f':\Sg \to N_n$ by pushing off $f(\mathcal{D}_f)$ to $\Int N_{n,\thick}$ such that the diameter of $f'(F)$ is bounded by a constant independent 
of $n$ for any component $F$ of $\ol{\Sg\setminus \mathcal{A}_f}$.
Note that $\mathcal{A}_{f\mathrm{:p}}={(f')}^{-1}(N_{n,\ppthick})$ is a subset of $\mathcal{A}_f=(f')^{-1}(N_{n,\thick})$ and 
$f'(F')$ is a surface in $N_{n,\ppthick}$ with $f'(\part F')\subset \part N_{n,\ppthick}$ for any component $F'$ of 
$\ol{\Sg\setminus \mathcal{A}_{f\mathrm{:p}}}$.
We say that the restriction $f'|\ol{\Sg\setminus \mathcal{A}_{f\mathrm{:p}}}$ is a \emph{pseudo-pleated map} and, for any component $F'$ of 
$\ol{\Sg\setminus \mathcal{A}_{f\mathrm{:p}}}$, $f'(F')$ is a 
\emph{pseudo-pleated surface}.
Since $f(\Sg)\subset C_n$, $f'(F')$ is contained in $C_{n,\ppthick}$.

\begin{lemma}\label{bdd_diam}
For any $R>0$, there exists a constant $d=d(R)>0$ independent of $n$ such that the diameter of any pseudo-pleated surface $H$ 
in $C_{n,\ppthick}$ with $H\cap \cb_R(x_n,N_n)\neq \eset$ is less than $d$.
\end{lemma}
\begin{proof}
Note that the number of components of $\Int H\cap N_{n,\mathrm{thin}}$ is at most $2(g-1)$, where $g$ is the genus of $\Sg$.
Moreover, the diameter of each component of $H\setminus N_{n,\mathrm{thin}}$ is bounded by a constant independent of $n$.
Since any component $V_n$ of $N_{n,\mathrm{thin}}$ meeting $\Int H$ non-trivially is not contained in $N_{n,\ppthin}$, 
$g_n|V_n:V_n\to M_{\mathrm{thick}}$ is an embedding onto a Margulis component of $M_{\mathrm{thin}}$ for all 
sufficiently large $n$.
There exist only finitely many Margulis tube components of $M_{\mathrm{thin}}$ which 
intersect non-trivially the $g_n$-image of at least one pseudo-pleated surface $H$ as above.
Hence we have a constant $c(R)>0$ independent of $n$ such that $\mathrm{diam}(V_n)<c(R)$.
This implies the existence of a constant $d(R)$ satisfying our desired condition.
\end{proof}

By an argument similar to that in the proof of \cite[Theorem 9.5.13]{tl}, for any point $y$ of $C_n$, there exists a pleated map $f:\Sg\to C_n\subset N_n$ with $f(\Sg)\cap \Int \cb_1(y,N_n)\neq \eset$.
The reason why one can use here the ``1-ball'' $\cb_1(\cdot)$ is explained in Soma \cite[p.\ 186]{so}.
In particular, for any $x\in C_{\pthick}$ and all sufficiently large $n$, we have a pseudo-pleated surface $F_n$ in $C_{n:\ppthick}$ meeting $\cb_1(g_n^{-1}(x),C_n)$ non-trivially.
By the Ascoli-Arzel\`{a} Theorem and Lemma \ref{bdd_diam}, $\{g_n(F_n)\}$ has a subsequence converging geometrically to a proper surface $F$ in $C_{\pthick}$ with $\mathrm{dist}(F,x)\leq 1$.
Such a surface is called a \emph{pseudo-pleated surface} in $C_{\pthick}$.
Suppose that $y$ is a point of some component $Q_n$ of the parabolic locus $P_n= P_n(\delta)$.
Then there exists a point $z$ of the corresponding component $Q_n(\ve)$ of $P_n(\ve)$ such that $\mathrm{dist}(y,z)$ is less than a constant depending only on $\delta$ (and $\ve=\ve(\delta)$).
Let $f:\Sg\to C_n$ be a pleated map with $\mathrm{dist}(f(\Sg),z)\leq 1$.
Since $f(\mathcal{D}_f)\cap \cn_1(N_{n,(p):\mathrm{thin}(\ve)})=\eset$, we have a pseudo-pleated map $f'|:\ol{\Sg\setminus \mathcal{A}_{f\mathrm{:p}}}\to C_{n,\ppthick}$ such that, for some component 
$F'$ of $\ol{\Sg\setminus \mathcal{A}_{f\mathrm{:p}}}$, the pseudo-pleated surface $f'(F')$ contains a point $w$ of $Q_n$ 
with $\mathrm{dist}(w,y)\leq d(\delta)$, where $d(\delta)>0$ is a constant depending only on $\delta$.
It follows that, for any point $x$ of any component $Q$ of $P=P(\delta)$, there exists a pseudo-pleated surface in 
$C_{\pthick}$ containing a point $v\in Q$ with $\mathrm{dist}(v,x)\leq d(\delta)$.

\subsection{Construction of block decomposition structures}

We often replace embedded surfaces by pseudo-pleated surfaces even though they are not necessarily embeddings.
Almost  families of pseudo-pleated surfaces considered in this section are equicontinuous  and have subsequences converging to 
other pseudo-pleated surfaces.
In order to proceed our argument, we need to restore an embedded surface from the limit pseudo-pleated surface.
Our procedure is represented as follows.
\begin{align*}
&\boxed{\mbox{Embedded surfaces}}\to\boxed{\mbox{pseudo-pleated surfaces}}\\
&\to\boxed{\mbox{a limit pseudo-pleated surface}}\to \boxed{\mbox{an embedded surface}}
\end{align*}

\begin{figure}[hbtp]
\centering
\scalebox{0.5}{\includegraphics[clip]{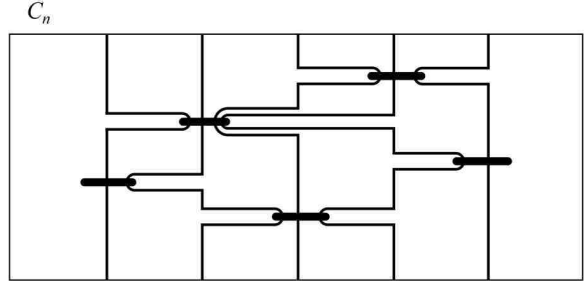}}
\caption{The union of bold horizontal segments is $C_{n,\ppthin}$.}
\label{H_n}
\end{figure}

Let $\ch_n$ be a finite set of mutually disjoint surfaces in $\mathrm{Int}C_n$ homeomorphic to $\Sg$ and such that each 
component of $\bigcup \ch_n\setminus \mathrm{Int}C_{n,\ppthick}$ is a properly embedded surface in 
$C_{n,\ppthick}$ and each component of $\Int C_{n,\ppthick}\setminus \bigcup \ch_n$ is homeomorphic to $\Sg\times (0,1)$, 
where $\bigcup \ch_n$ denotes the union $\bigcup_{H_\alpha\in \ch_n}H_\alpha$ as usual.
Let $\ch_n^{(1)}$ be the set of components of $\bigcup \ch_n\setminus C_{n,\ppthick}$.
If $H\in \ch_n^{(1)}$ has a non-peripheral and non-contractible (i.e.\ \emph{essential}) simple loop $l$ properly isotopic in 
$C_{n,\ppthick}$ to a loop $l'$ in the parabolic locus $P_n$.
By performing surgery on $H$ along an embedded annulus $A$ in $C_{n,\ppthick}$ with $\part A=l\cup l'$, we have a properly embedded 
surface $H'$ in $C_{n,\ppthick}$ homeomorphic to $H\setminus \Int \cn(l)$ and $\part H'\setminus \part H \subset P_n$, where $\cn(l)$ is a regular neighborhood of $l$ in $H$.
Repeating the same argument repeatedly, we have a finite set $\ch_n^{(2)}$ of mutually disjoint, properly embedded 
surfaces in the pared manifold $(C_{n,\ppthick},P_n)$ satisfying the following conditions.
\begin{enumerate}[\rm (i)]
\item
For each element $H'$ of $\ch_n^{(2)}$, any essential simple loop in $H'$ is not homotopic in $C_{n,\ppthick}$ 
to a loop in $P_n$.
\item
Each component of $\Int C_{n,\ppthick}\setminus \bigcup\ch_n^{(2)}$ is a block, see Fig.\ \ref{H_n}.
\end{enumerate}
We may assume that any two elements of $\ch_n^{(2)}$ are not mutually parallel in $C_{n,\ppthick}$ 
if necessary by taking a representative in each mutually parallel class in $\ch_n^{(2)}$. 

Each $H_\alpha'\in \ch_n^{(2)}$ is properly homotopic in $C_{n,\ppthick}$ to a pseudo-pleated surface $F_\alpha$.
The set $\{F_\alpha\}_{H_\alpha'\in \ch_n^{(2)}}$ is denoted by $\cf_n$.
Here, we note that, in general, one can not detect such a surface $F_\alpha$ with a uniformly bounded distance from $H_\alpha'$.
In fact, though the ubiquity of pseudo-pleated surfaces guarantees the existence of a pseudo-pleated surface $F_\alpha'$ near $H_\alpha'$, 
it does not necessarily imply that $F_\alpha'$ is properly homotopic to $H_\alpha'$ in $C_{n,\ppthick}$.
Set $C_{\pthick}^\circ=C_{\pthick}\setminus P$ and $C_{n,\ppthick}^\circ=C_{n,\ppthick}\setminus P_n$.
The intersection $F_\alpha\cap C_{n,\ppthick}^\circ$ is denoted by $F_\alpha^\circ$, and the family $\{F_\alpha^\circ; F_\alpha\in \cf_n\}$ by $\cf_n^\circ$.

When $J$ is a proper surface in a 3-manifold properly homotopic to an embedded surface,  
we denote the surface by  $J^\flat$ and the family $\{J_\alpha^\flat\}_{J_\alpha\in \mathcal{J}}$ by $\mathcal{J}^\flat$.

\begin{lemma}\label{embeddable}
Let $J_n$ $(n=1,2,\dots)$ be pseudo-pleated surfaces in $C_{n,\ppthick}$ each of which is properly homotopic to an embedded surface and such that $\{g_n(J_n)\}$ converges to 
a pseudo-pleated surface $J$ in $C_{\pthick}$.
Then $J$ is properly homotopic in $C_{\pthick}$ to an embedded surface $J^\flat$. 
\end{lemma}
\begin{proof}
By some arguments in Freedman, Hass and Scott \cite{fhs}, each $J_n$ is properly homotopic to an embedding $J_n^\flat$ in $C_{n,\ppthick}$ contained in an arbitrarily small neighborhood of $J_n$ and the support of the homotopy is in a small neighborhood $Z_n$ of the union of $J_n\cup J_n^\flat$ and all bounded components of $C_{n,\ppthick}\setminus (J_n\cup J_n^\flat)$.
From the ubiquity of pseudo-pleated surfaces, for any $x\in Z_n\setminus \cn_1(J_n\cup J_n^\flat)$, there exists a 
pseudo-pleated surface in $C_{n,\ppthick}$ meeting both $\cb_1(x,N_n)$ and $J_n\cup J_n^\flat$ non-trivially.
By Lemma \ref{bdd_diam}, the diameter of $Z_n$ is bounded by a constant independent of $n$.
Fix $n\in \nn$ with $Z_n\subset \cb_{R_n}(x_n,N_n)$.
Then $J$ is properly homotopic to an embedded surface $J^\flat$ in $g_n(Z_n)\subset C_{\pthick}$.  
\end{proof}

The following proposition is a key to the proof of Theorem \ref{main_a}.

\begin{prop}\label{M_block}
$M^\natural$ admits a block decomposition.
\end{prop}
\begin{proof}
By Lemma \ref{bdd_diam} and the Ascoli-Arzel\`{a} Theorem, there exists a (possibly empty) set $\cf$ of pseudo-pleated 
surfaces in $C_{\pthick}$ such that $g_n(\cb_{R_n}(x_n,N_n)\cap \bigcup \cf_n)=\cb_{R_n}(x_\infty,
M)\cap \bigcup \cf$ if necessary passing to a subsequence of $\{g_n\}$ and replacing $\{R_n^{\mathrm{old}}\}$ by 
some divergence sequence $\{R_n^{\mathrm{new}}\}$ with $0<R_n^{\mathrm{new}}<R_n^{\mathrm{old}}$.
By Lemma \ref{embeddable} together with the least area surface theory in \cite{fhs}, for any $F_\alpha\in \cf$,  
there exists an embedded surface $F_\alpha^\flat$ in $C_{\pthick}$ such that $F_\alpha^\flat\cap F_\beta^\flat=\eset$ for any two distinct 
elements $F_\alpha,F_\beta\in \cf$.

Since $C_{\pthick}^\circ$ is homeomorphic to $C$ and hence to $M^\natural$, it suffices to show that    
$C_{\pthick}^\circ$ is divided into blocks along locally finite open g-subsurfaces.
Since any finite union of such blocks is realized as a subspace of $C_n\cong \Sg\times I$ via $g_n^{-1}$ with sufficiently large $n$, by Proposition \ref{key} there exists an embedding $h:M^\natural\cong C_{\pthick}^\circ\to \Sg\times I$ 
whose image is the support of a block complex.

First, we will prove that any component $B$ of $C_{\pthick}^\circ\setminus \bigcup{\mathcal F}^{\flat\circ}$ is the support of a finite block complex.
Let $\cf^{\flat\circ}_{n,B}$ be the subset of $\cf_n^{\flat\circ}$ consisting of $F_{\alpha}^{\flat\circ}\in \cf_n^{\flat\circ}$ with $F_{\alpha}^{\flat\circ}\subset \cb_{R_n}(x_n,N_n)$ and 
$g_n(F_{\alpha}^{\flat\circ})\subset \mathrm{Fr}(B)$.
Since $\cf_{n,B}^{\flat\circ}$ consists of mutually non-parallel surfaces contained in the frontier of a single block, the number of elements of 
$\cf_{n,B}^{\flat\circ}$ is bounded by a constant independent of $n$.
Thus, for all sufficiently large $n\in\nn$, the images $g_n(\bigcup \cf_{n,B}^{\flat\circ})$ are the same 
subset of $\cf^{\flat\circ}$.
Since each component of $C_{n,\ppthick}^\circ\setminus \bigcup \cf_n^{\flat\circ}$ is a block, there exists a subset 
$\cu_{n,B}^\circ$ of $\cf_n^\circ\setminus \cf_{n,B}^\circ$ such that $\bigcup (\cf_{n,B}^{\flat\circ} \cup \cu_{n,B}^{\flat\circ})$ excises from $C_{n,\ppthick}^\circ$ a 
block $B_n$ with $\mathrm{Fr}(B_n)=\bigcup (\cf_{n,B}^{\flat\circ} \cup \cu_{n,B}^{\flat\circ})$, where we choose $\cu_{n,B}^{\flat\circ}$ so that $\bigcup \cu_{n,B}^{\flat\circ}$ is contained in an arbitrarily small neighborhood of $\bigcup \cu_{n,B}^\circ$ in $C_{n,\ppthick}^\circ$.
Since $\cu_{n,B}^\circ$ is not contained in $\cb_{R_n}(x_n,N_n)$, for any $R>0$, there exists $n_0=n_0(R)\in \nn$ such that 
$\cn_R(S_n)\cap \cu_{n,B}^{\flat\circ}=\emptyset$ if $n\geq n_0$, where $S_n=\bigcup \cf_{n,B}^{\flat\circ}$ if $\cf_{n,B}^{\flat\circ}\neq \eset$ otherwise $S_n$ is the one-point set $\{g_n^{-1}(y_\infty)\}$ for a fixed $y_\infty \in B$.
From the ubiquity of pseudo-pleated surfaces, there exists a finite set ${\mathcal{W}}_{n,B}^\circ$ of pseudo-pleated surfaces in $B_n$ such that $B_n\setminus \bigcup {\mathcal{W}}_{n,B}^{\flat\circ}$ contains a block component $B_n'$ 
as a deformation retract of $B_n$ and satisfying 
$$r_1<\mathrm{dist}\Bigl(\bigcup {\mathcal{W}}_{n,B}^\circ,S_n \Bigr)<r_2\quad 
\mbox{and hence}\quad\lim_{n\sto \infty}\mathrm{dist}\Bigl(\bigcup {\mathcal{W}}_{n,B}^\circ,\bigcup {\mathcal{U}}_{n,B}^{\flat\circ}\Bigr)=\infty$$
for some constants $0<r_1<r_2$ independent of $n$, see Fig.\ \ref{s_n}. 
\begin{figure}[hbtp]
\centering
\scalebox{0.5}{\includegraphics[clip]{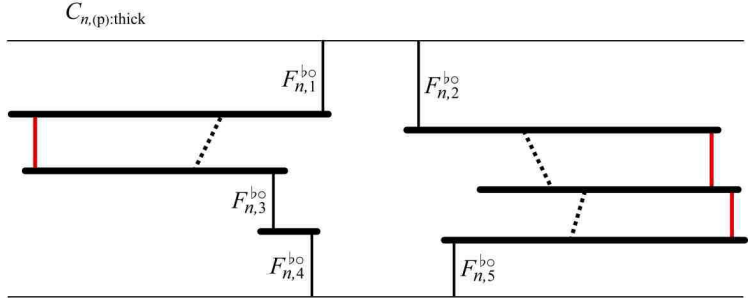}}
\caption{The union of vertical red segments is $\bigcup {\mathcal{U}}_{n,B}^{\flat\circ}$ and 
that of dotted segments is $\bigcup {\mathcal{W}}_{n,B}^\circ$.}
\label{s_n}
\end{figure}
In particular, one can obtain a finite set $\mathcal{W}_B$ of pseudo-pleated surfaces in $C_{\pthick}^\circ$ 
such that $\{g_n(\bigcup {\mathcal{W}}_{n,B}^\circ)\}$ has a subsequence of converging uniformly to $\bigcup {\mathcal{W}}_B^\circ$.
By Lemma \ref{embeddable}, there exists a family $\mathcal{W}_B^{\flat\circ}$ of properly embedded surfaces in 
$C_{\pthick}^\circ$.
Then $\bigcup {\mathcal{W}}_B^{\flat\circ}$ excises from $B$ components adjacent to geometrically infinite 
ends.
In particular, the union $E$ of these components is homeomorphic to 
$\bigcup {\mathcal{W}}_B^{\flat\circ}\times [0,\infty)$.
We may choose $\mathcal{W}_{n,B}^{\flat\circ}$ so that $g_n(\bigcup \mathcal{W}_{n,B}^{\flat\circ})=\bigcup \mathcal{W}_B^{\flat\circ}$ 
for all sufficiently large $n$.
Since the diameter of $\fr(B_n')=\bigcup( \cf_{n,B}^{\flat\circ}\cup \mathcal{W}_{n,B}^{\flat\circ})$ is uniformly bounded, as in the 
proof of Lemma \ref{embeddable}, one can show that the diameter of $B_n'$ is also uniformly bounded.
Fix $n\in \nn$ with $B_n'\subset \cn_{R_n}(x_n,N_n)$.
Then the complement $B\setminus E=g_n(B_n')$ is also a block in $C_{\pthick}^\circ$.

Now,  we know that the union of ${\mathcal F}^{\flat\circ}$ and $\mathcal{W}_B^{\flat\circ}$'s 
defines a block decomposition $\mathcal{L}$ on $C_{\pthick}^\circ$.
Let $\hat B_1,\dots,\hat B_k$ be the elements of $\mathcal{L}$ meeting $F_0$ \emph{essentially}, that is, 
any surface in $C_{\pthick}^\circ$ homotopic to $F_0$ intersects $\hat B_i$ non-trivially.
Our desired block complex $\ck$ is obtained from $\mathcal{L}$ by replacing $\{\hat B_1,\dots,\hat B_k\}$ by a single block $B_0=\Int (\ol{\hat B_1\cup\cdots\cup \hat B_k})$, which is considered to be the base block of $\ck$.
\end{proof}

\section{Proof of Theorem \ref{main_a}}\label{S_4}

We will still work with the notation and definitions in Section \ref{S_3}.
By Proposition \ref{M_block}, there 
exists an embedding $\psi:M^\natural=\mathbf{H}^3\cup \O(G)/G\to \Sg\times I$ with 
$\psi^{-1}(\Sg_0\cup \Sg_1)=\part M^\natural$ and such that 
each $X_y=\cx\cap \Sg_y$ $(y \in \cy=q(\cx))$ is a g-l-subsurfaces of $\Sg_y$, where $\cx=\Sg\times I\setminus \psi(M^\natural)$.
Thus $\cx$ satisfies the condition (i) of Theorem \ref{main_a}.
If necessary saturating $\cx$, we may assume that the set $\cp_\cx$ defined as in Subsection \ref{ss_2_3} 
consists of horizontal annuli in $\Sg\times I$.
From our construction of the block complex on $M^\natural$, any element of $\cp_\cx$ corresponds to a parabolic cusp of $M$.
Throughout this section, the image $\psi(M^\natural)=\Sg\times I\setminus \cx$ is denoted by $L$.

\begin{lemma}\label{condition_ii}
$\cx$ satisfies the condition {\rm (ii)}.
\end{lemma}
\begin{proof}
Let $F$ be any non-peripheral front of $\cx$.
If $F$ were homeomorphic to an open pair of pants, then a pleated map $F\to M$ corresponding to the inclusion of $F$ to the closure $\ol L$ of $L$ in $\Sg\times I\setminus \cp_\cx$ 
would be a totally geodesic embedding the image of which excises from $M$ a submanifold $E$ facing a geometrically finite end.
This means that $E$ is adjacent to a component of $\part M^\natural$, which contradicts that 
$\psi^{-1}(\Sg_0\cup \Sg_1)=\part M^\natural$.
\end{proof}

\begin{proof}[Proof of Theorem \ref{main_a}]
To complete the proof, it remains to show that $\cx$ satisfies the condition (iii).

For the proof, we use the fact that any parabolic cusps of $M$ are not parallel to each other.
Since each component of $\cp_\cx$ corresponds to a parabolic cusp of $M$, it suffices to show that any component 
$l_y$ of $\lambda(\L_y^\varepsilon)\setminus \fr(X_y)$ $(\varepsilon=\pm)$ corresponds to 
a parabolic cusp of $M$.
If $l_y\cap \L_y^+\neq \eset$, then by Lemma \ref{l21} (and Remark \ref{R2}), $l_y$ is contained in an element of $\cp_\cx$, and 
hence $M$ has a parabolic cusp corresponding to $l_y$.
Thus we may assume that $l_y\cap \L^+_y=\eset$.
Then there exists a sequence $\{H_n\}$ of vertical open g-subsurface in $L$ with $\fr(H_n)\subset \cx$ approaching to $l_y$ and 
such that any two $H_n,H_m$ of them are not properly homotopic to each other in $L$, see Fig.\ \ref{hidden_p}\,(a).
\begin{figure}[hbtp]
\centering
\scalebox{0.5}{\includegraphics[clip]{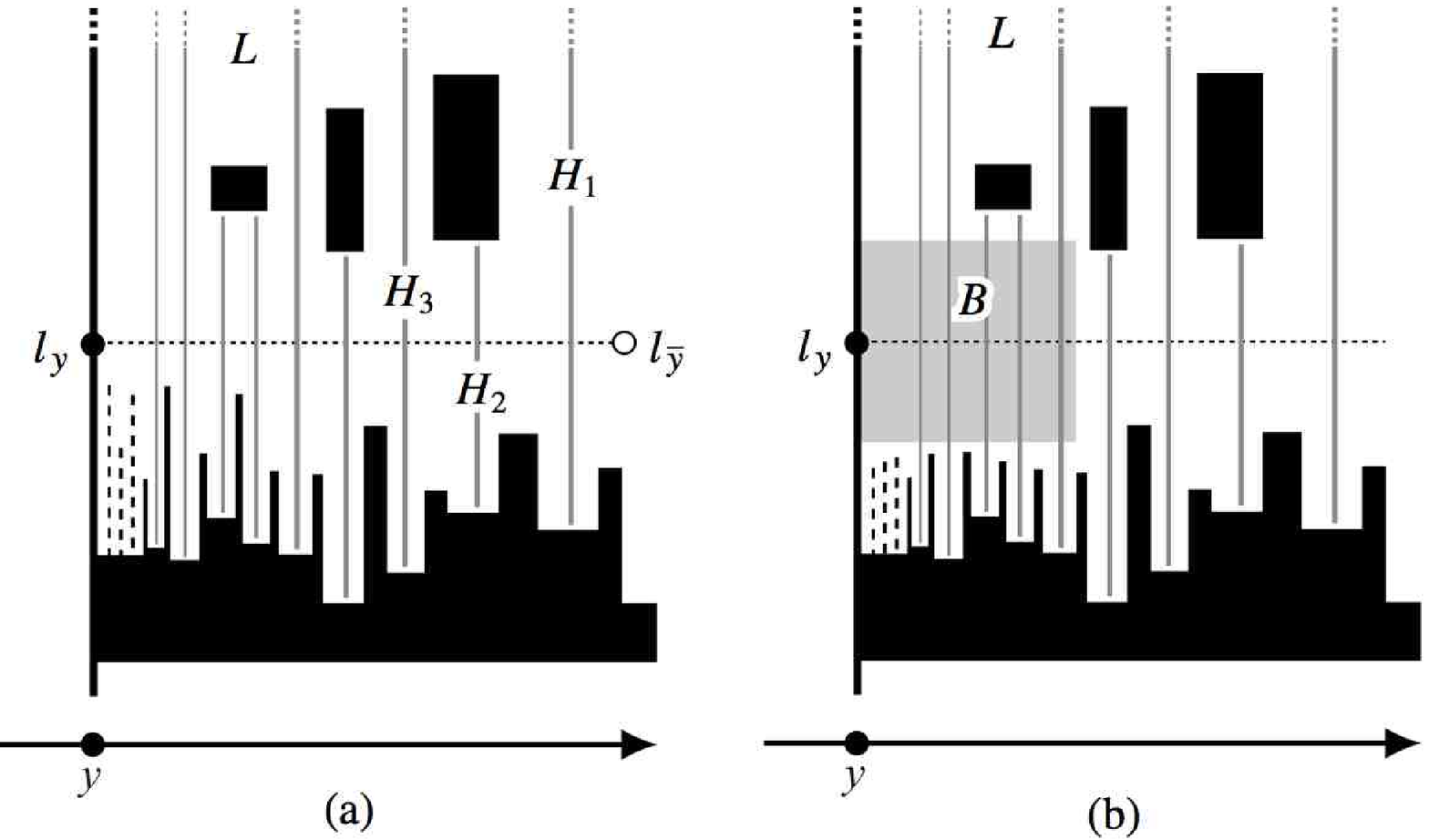}}
\caption{}
\label{hidden_p}
\end{figure}

We suppose that $M$ does not have a parabolic cusp corresponding to $l_y$ and introduce a contradiction.

First consider the case when there exists a geodesic loop $\lambda$ in $M$ corresponding to a g-loop $l_{\bar y}$ in $L\cap \Sg\times (y,1]$ 
parallel to $l_y$ in $L$.
If necessary retaking $\delta>0$, we may assume that $\lambda\cap M_{\pthin}=\eset$.
Let $f_n:H_n\to M$ be a pleated surface corresponding to the inclusion $H_n\subset L$ and with 
$f_n(H_n)\supset \lambda$.
By the Ascoli-Arzel\`{a} Theorem, there exists an open  g-subsurface $J$ of $\Sg$ and an embedding $i_n:J\to H_n$ such 
that $J_n=i_n(J)$ is an open g-subsurface of $H_n$ containing the g-loop $l_{y,n}:=(l_{y})_{\{y_n\}}$ in $H_n$ and $\{f_n\circ i_n\}$ is a 
sequence of proper maps to $(M,M_{\mathrm{thin}})$ converging  uniformly a continuous proper map if necessary passing to a subsequence. 
In particular, there exists $n_0\in  \nn$, such that $f_n\circ i_n$ $(n\geq n_0)$ are homotopic to each other in $M$.
This implies that the open g-subsurfaces $J_n$ of $H_n$ are horizontally parallel to each other in $L$.
It follows that $L$ contains a block $B$ of the form $\hat J_{(y,y+\varepsilon)}$ for a sufficiently small $\varepsilon>0$ (see Fig.\ \ref{hidden_p}\,(b)), 
where $\hat J=\mathrm{pr}_\mathrm{hz}(i_n(J_n))$ is an open g-subsurface of $\Sg$ containing $\mathrm{pr}_\mathrm{hz}(l_y)$ and independent of $n\geq n_0$.
This contradicts that $l_y$ is a component of $\lambda(\L_y^+)\setminus \fr(X_y)$.

Next consider the case when there exists a component $P$ of $M_{\pthin}$ corresponding to a horizontal annulus or a g-loop in 
$\cx\cap \Sg\times (y,1]$ parallel to $l_y$ in $L$.
Then one can again have a contradiction by considering pleated maps $f_n':H_n\setminus l_{y,n}\to M$ instead of $f_n$ as above 
such that small neighborhoods of the both ends of $H_n\setminus l_{y,n}$ adjacent to $l_{y,n}$ are mapped into $P$.
This completes the proof of Theorem \ref{main_a}.
\end{proof}

\section{Proof of Theorem \ref{main_b}}\label{S_5}

In Subsection \ref{ss_5_1}, we give the definition of efficient pleated surfaces.
Subsection \ref{ss_5_3} presents several sequences of Kleinian groups 
which give geometric limits $G$ and $G'$ in two manners.
In Subsection \ref{ss_5_4}, we show that (i) $\mathbf{H}^3\cup \O(G)/G$ have a topological type desired in Theorem \ref{main_b}, (ii) $G=G'$  and (iii) $G'$ is a geometric limit of quasi-Fuchsian groups.

We say that a sequence $\{a_n\}$ in a metric space $A$ \emph{subconverges} to $a\in A$ if $\{a_n\}$ has a subsequence 
converging to $a$.
For simplicity, such a subsequence is denoted again by $\{a_n\}$.

\subsection{Length of measured laminations and efficient pleated maps}\label{ss_5_1}

Let $F$ be either $\Sg$ itself or an open g-subsurface of $\Sg$.
For any $\mu \in \mathcal{ML}_0(F)$ and $\sigma\in \mathrm{Teich}(F)$, we have g-loops $l_n$ in $F(\sigma)$ and $t_n>0$ such that the sequence $\{t_nl_n\}$ of weighted loops converges to $\mu$ in $\mathcal{ML}_0(F)$.
Then there exists the limit $\mathrm{length}_{F(\sigma)}(\mu)=\lim_{n\sto \infty}t_n\mathrm{length}_{F(\sigma)}(l_n)$, called the {\it length} of $\mu$ in $F(\sigma)$, which is independent of the choice of a sequence $\{t_nl_n\}$ converging to $\mu$.
According to Bonahon \cite[Proposition 4.5]{bo}, for $\mu=\lim_{n\sto \infty}t_nl_n,\mu'=\lim_{n\sto \infty}t_n'l_n'\in \mathcal{ML}_0(F)$, the {\it geometric intersection number} $i(\mu,\mu')$ is well defined by $\lim_{n\sto \infty} t_nt_n' i(l_n,l_n')$, where $i(l_n,l_n')$ is the cardinality of $l_n\cap l_n'$ when they meet each other transversely and otherwise $i(l_n,l_n')=0$.

\begin{lemma}\label{l5_1}
Suppose that $F$ has a completed hyperbolic structure $\sigma$ with finite area and $H$ is a 
proper open g-subsurface of $F$. 
Let $\nu$ be a measured lamination in $H$.
Then, for any connected measured lamination $\mu$ in $F$ with $i(\mu,\nu)\neq 0$, there exists a measured lamination $\mu'$ in $H$ independent of $\sigma\in \mathrm{Teich} (F)$ with 
$i(\mu',\nu)\neq 0$ and 
$$\mathrm{length}_H(\mu')\leq i(\mu,\mathrm{Fr}(H))\mathrm{length}_{F(\sigma)}(\mathrm{Fr}(H))+
2\mathrm{length}_{F(\sigma)}(\mu),$$
where $H$ are supposed to have the incomplete hyperbolic metric induced from $F(\sigma)$.
\end{lemma}
\begin{proof}
If $\mu\subset H$, then we may set $\mu'=\mu$.
If not, $\mu$ meets $\mathrm{Fr}(H)$ transversely and non-trivially.
Since $\mu$ is connected, each leaf of $\mu$ intersects $\mathrm{Fr}(H)$.
In particular, a leaf $s$ in $\mu|H$ with $s\cap \nu\neq \eset$ is an open geodesic arc with $\ol s\cap \mathrm{Fr}(H)\neq \eset$, 
where $\ol s$ is the closure of $s$ in $F$.
Let $c$ be the union of components of $\mathrm{Fr}(H)$ meeting $\ol s$ non-trivially.
The $c$ consists of one or two g-loops in $F$.
Let $\alpha$ be the union of all leaves in $\mu |H$ properly homotopic to $s$ in $H$, and let $s_0$ be the shortest geodesic arc in the proper homotopy class of $s$.
The closure $\ol s_0$ of $s$ meets $c$ orthogonally in two points.
If we set $t=i(\ol\alpha,c)$, then $t\leq i(\mu,c)$ and $\mathrm{length}_H(ts_0)\leq 
2\mathrm{length}_H(\alpha)\leq 2\mathrm{length}_F(\mu)$.
Since $H$ contains the measured lamination $\nu$, $H$ is not an open pair of pants.
It follows that, for a small regular neighborhood $N$ of $s_0 \cup c$ in $F\cup c$, $\partial N\cap H$ contains at least one non-peripheral loop $b$ in $H$.
Let $\lambda$ be the g-loops in $H$ homotopic to $b$.
The $\lambda$ is supposed to be a measured lamination with the counting measure.
Then one can choose the $b$ so that $i(\lambda, \nu)\geq i(s_0,\nu)$, see Fig.\ \ref{f_h}.
\begin{figure}[hbtp]
\centering
\scalebox{0.5}{\includegraphics[clip]{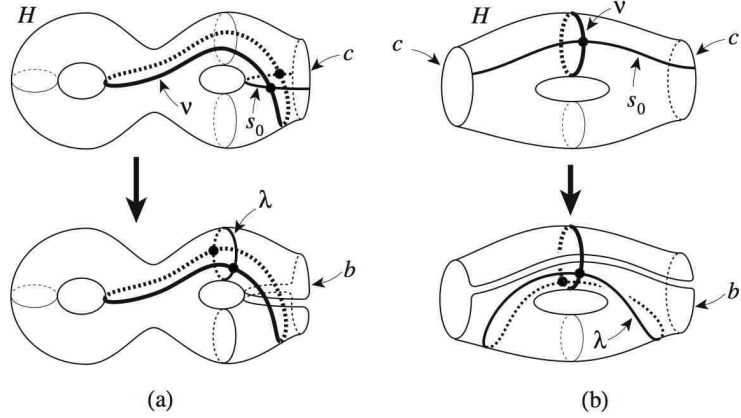}}
\caption{The case when $\nu$ is a closed geodesic.
In (a), $i(\lambda,\nu)=i(s_0,\nu)=2$.
In (b), $i(\lambda,\nu)=2$, $i(s_0,\nu)=1$.}
\label{f_h}
\end{figure}
If we set $\mu'=(t/2)\lambda$, then $i(\mu',\nu)=(t/2) i(\lambda,\nu)\geq (t/2)i(s_0,\nu)\neq 0$.
Moreover,
\begin{align*}
\mathrm{length}_F(\mu')&=\frac{t}{2}\mathrm{length}_F(\lambda)\\
&<\frac{t}{2}\cdot 2(\mathrm{length}_{F}(c)+\mathrm{length}_H(s_0))\\
&\leq i(\mu,c)\mathrm{length}_{F}(c)+2\mathrm{length}_F(\mu).
\end{align*}
This implies our desired inequality.
\end{proof}

Let $M$ be a hyperbolic 3-manifold admitting a proper homotopy equivalence $h:F\to M$, called a \emph{marking}.
The marking $h$ induces the holonomy $\rho:\fd(F)\to \mathrm{PSL}_2(\cc)$.
If a connected measured lamination $\mu$ in $F$ is not realizable in $M$ by any pleated map properly homotopic to $h$, then we set $\mathrm{length}_\rho(\mu)=0$.
If there exists a pleated map $f_\mu:F(\sigma)\to M$ properly homotopic to $h$ and realizing $\mu$, then we set $\mathrm{length}_\rho(\mu)=\mathrm{length}_{F(\sigma)} (\mu)$.
Since $f_\mu:\Sg(\sigma)\to M$ is arcwise isometric, $\mathrm{length}_\rho (\mu)$ represents the length of a measured lamination $\ol\mu=f_\mu(\mu)$ in $M$, that is, $\mathrm{length}_\rho (\mu)=\mathrm{length}_M (\ol\mu)$.
If $\nu\in \mathcal{ML}_0(F)$ is decomposed to the sum of connected sublaminations $\nu =\mu_1 +\cdots +\mu_n$, then we set $\mathrm{length}_\rho(\nu)=\mathrm{length}_\rho (\mu_1)+\cdots +\mathrm{length}_\rho (\mu_n)$.

Consider a sequence $\{\rho_n:\fd(F)\to \mathrm{PSL}_2(\cc)\}$ of faithful discrete representations transferring any peripheral element of $\fd(F)$ to a parabolic element.
In the case when $F\neq \Sg$, fix a complete hyperbolic structure $\sigma$ on $F$ of finite area.
Let $\lambda$ be a disjoint union of finitely many simple geodesic lines in $F(\sigma)$ such that each component of $F(\sigma)\setminus\lambda$ is an open ideal triangle.
The $\lambda$ is regarded as a `lamination' in $F(\sigma)$ with non-compact support.
In the case when $F=\Sg$, consider a g-loop $l$ in $\Sg$ with $d=\inf_n \{\mathrm{length}_{\rho_n}(l)\}>0$, see \cite[Proposition 4.2]{tii}.
Let $\lambda'$ be a disjoint union of finitely many simple geodesics lines in $\Sg \setminus l$ such that each component of $\Sg \setminus l\cup \lambda'$ is an open ideal triangle.
Then $\lambda=l\cup \lambda'$ is a geodesic lamination in $\Sg$ consisting of a single compact leaf $l$ and finitely many non-compact leaves spiraling around $l$.
Let $f_{\lambda,n}:F(\sigma_n)\to N_n=\mathbf{H}^3/\rho_n(\fd(F))$ be a pleated map realizing $\lambda$.
According to \cite[Theorem 3.3]{tii}, there exists a constant $C=C(d)>0$ depending only on $d$ and a continuous non-negative function $a_\lambda :\mathcal{ML}_0(F)\to \rr$ depending only on $\lambda$ such that, for any $\mu\in \mathcal{ML}_0(F)$ and $n\in \nn$,
\begin{equation}\label{length}
\mathrm{length}_{\rho_n}(\mu)\leq \mathrm{length}_{F(\sigma_n)}(\mu)\leq \mathrm{length}_{\rho_n}(\mu)+Ca_\lambda (\mu).
\end{equation}
This inequality is very useful for examining whether a given sequence $\{\rho_n\}$ has an algebraically convergent subsequence.
We say that these $f_{\lambda,n}$ are \emph{efficient pleated maps} for $\{\rho_n\}$.

\subsection{Construction of approaching Kleinian groups}\label{ss_5_3}

Throughout the remainder of this section, we fix a crevasse $\cx$ satisfying the conditions same to those on $\hat \cx$ of Lemma \ref{c_to_b} 
and set $\mathcal{S}=\Sg_0\cup \Sg_1\setminus \cx$.
Moreover, as in Subsection \ref{p_c}, we may assume that any non-peripheral front of $\cx$ is not in $\Sg_0\cup \Sg_1$.
Let $\cp_\cx$ be the set of horizontal annuli in $\cx$ defined in Subsection \ref{ss_2_3}.
As is suggested in Fig.\ \ref{x_hat1}\,(b) and \ref{x_hat2}\,(b), for any component $l$ of $\lambda(\L_y^\ve)$ $(y\in \cy, \ve=\pm)$, 
there exists an element of $\cp_\cx$ containing $l$ as an edge.

Any Kleinian 3-manifold $N^\natural =\mathbf{H}^3\cup \O(\G)/\G$ treated in this section has a marking from a subset $L'$ of $\Sg\times I$ with 
$L'\supset \Sg_{1/2}\cup \mathcal{S}$ and 
such that a continuous map $f_N:\Sg\to N$ corresponding the the inclusion $\Sg_{1/2}\subset L'$ is $\fd$-injective and 
there exists the union $\part_\cs N^\natural$ of components of $\part N^\natural$ corresponding to $\mathcal{S}$.
Moreover, there exists a full measured lamination $\mu$ on $\Sg$  without compact leaves which is taken commonly for all $N$ 
and such that 
$f_N$ is homotopic to a pleated surface $f_{N,\mu}$ realizing $\mu$ in $N$.
We always take the base point $x$ of $N$ in $\mathrm{Image}(f_{N,\mu})$.
If $p:N_1=\mathbf{H}^3/\G_1\to N=\mathbf{H}^3/\G$ is a locally isometric covering between hyperbolic 3-manifolds such that $p\circ f_{N_1}$ is 
homotopic to $f_{N}$, then $p\circ f_{N_1,\mu}=f_{N,\mu}$.
Thus one can choose the base point $x_1$ of $N_1$ so that $p(x_1)=x$.
The covering $p$ is naturally extended to that over $N^\natural$.
Since $\O(\G)\subset \O(\G_1)$, the total space $\mathbf{H}^3\cup \O(\G)/\G_1$ of $p^\natural$ is a subset of $N_1^\natural=\mathbf{H}^3\cup \O(\G_1)/\G_1$.
In general, its boundary $\O(\G)/\G_1$ is a proper subset of $\part N^\natural=\O(\G_1)/\G_1$.
However, for any covering $p$ considered below, $\O(\G)/\G_1$ contains $\part_\cs N_1^\natural$ 
as a union of components and 
the restriction  $p^\natural |\part_\cs N_1^\natural$ is a marking-preserving homeomorphism to $\part_\cs N^\natural$.
The following diagram represents the situation.
$$
\begin{CD}
\part_\cs N_1^\natural @>\text{inclusion}>> \O(\G)/\G_1@>\text{inclusion}>>\part N_1^\natural\\
@VV{\cong}V @VV{p^\natural|\O(\G)/\G_1}V\\
\part_S N^\natural @>\text{inclusion}>> \O(\G)/\G
\end{CD}
$$

By Lemma \ref{c_to_b}, $L=\Sg\times I\setminus \cx$ admits a block decomposition $\ck_\infty=\bigcup_{n=1}^\infty \ck_n$ such that any maximal horizontal open annulus in $L$ meets only finitely many blocks in $\ck_\infty$, 
where $\ck_n$'s are finite block complexes with $\ck_n\subset \ck_{n+1}$, possibly $\ck_n=\ck_{n+1}=\cdots$ for some  
$n\in \nn$ if $\fd(L)$ is finitely generated. 
One can also suppose that all blocks of $\ck_\infty$ containing components of $\cs$ belong to $\ck_1$.
The support of $\ck_n$ is denoted by $L_n$.
Then the union $\bigcup_{n=1}^\infty L_n$ is equal to $L$.

\begin{figure}[hbtp]
\centering
\scalebox{0.5}{\includegraphics[clip]{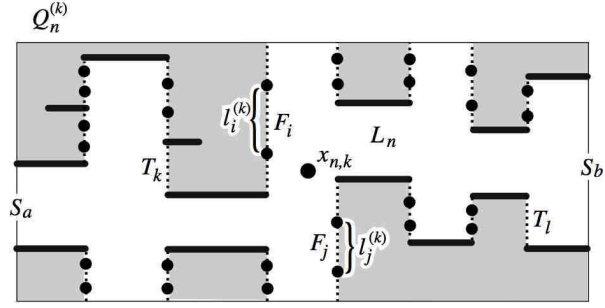}}
\caption{Each bold horizontal segment represents an element of $\cp_n$, and 
each vertical dotted segment with two bold dots does either a non-peripheral front of $\cx$ or an element of $\mathcal{J}_n\setminus \mathcal{T}_n$.
$T_k,T_l$ are elements of $\mathcal{T}_n$.
$S_a,S_b$ are peripheral fronts of $\cx$.}
\label{add_1}
\end{figure}

Let $\cp_n$ be the subset of $\cp_\cx$ consisting of elements meeting the closure $\ol L_n$ in $\Sg\times I$ nontrivially and let $\cz_n$ be the set of components of $\partial \ol L_n\setminus \bigcup\cp_n$.
By the condition (ii) of $\cx$, each element of $\cz_n$ not in $\cs$ is either a non-peripheral front of $\cx$ or a joint of a block in $\ck_n$ and that in $\ck_\infty\setminus \ck_n$.
Let $\mathcal{J}_n$ be the subset of $\cz_n$ consisting of such joints.
Each element of $\cz_n\setminus \mathcal{J}_n$ is not an open pair of pants, but $\mathcal{J}_n$ may contain them.
Let $\mathcal{T}_n$ be the subset of $\mathcal{J}_n$ consisting of elements homeomorphic to open pair of pants.
Thus we have
$$\mathcal{Z}_n\supset \mathcal{J}_n\supset \mathcal{T}_n.$$
For any non-peripheral element $F_j$ of $\cz_n\setminus \mathcal{T}_n$, take a full measured lamination $\mu_j$ on $F_j$ without compact leaves 
arbitrarily and a $k$-sequence $\{l_j^{(k)}\}$ of g-loops in $F_j$ such that the $k$-sequence $\{t_j^{(k)}l_j^{(k)}\}$ converges to $\mu_j$ in $\mathcal{ML}_0(F_j)$ for some $t_j^{(k)}>0$ $(k\in \nn)$.
For any peripheral front $S_a$ of $\cx$, take a conformal structure $\sigma_a$ on $S_a$ arbitrarily.
By the condition (iii) of $\cx$, it is not hard to show that  
$$Q_n^{(k)}=\Sg\times I \setminus \bigcup\bigl(\cp_n \cup \cp_n^{(k)}\bigr)$$
is atoroidal for all sufficiently large $k\in \nn$, 
where $\cp_n^{(k)}$ is the set of g-loops $l_j^{(k)}$ as above, see Fig.\ \ref{add_1}.
According to Thurston's Uniformization Theorem for atoroidal Haken 
manifolds (see Kapovich \cite{ka} for the complete proof), $Q_j^{(k)}$ admits a geometrically finite Kleinian structure such that each component of $\part Q_j^{(k)}$ 
corresponding to a peripheral element $S_a\in \cz_n$ has the conformal structure $\sigma_a$.
Moreover, the set of parabolic cusps of $Q_j^{(k)}$ corresponds to $\cp_n\cup \cp_n^{(k)}$ bijectively.
Let $\{F_1,\ldots,F_s\}$ be the set of non-peripheral elements of $\cz_n\setminus \mathcal{T}_n$.
By performing hyperbolic Dehn surgery on $Q_n^{(k)}$ along $l_1^{(k)},\dots,l_s^{(k)}$ of type $(1,u_k)$ with sufficiently large 
$u_k\in \nn$, one can have a geometrically finite Kleinian manifold $N_n^{(k)\natural}$ without changing the conformal structure on 
$\partial Q_n^{(k)}=\partial N_n^{(k)\natural}$ and such that the geodesic loops $\ol l_j^{(k)}$ in the hyperbolic 3-manifold $N_n^{(k)}$ 
corresponding to $l_j^{(k)}$ satisfy $\lim_{k\sto \infty}\mathrm{length}(\ol l_j^{(k)})=0$. 
Since $L_n$ is a subset of $\Sg\times I$ disjoint from $\bigcup \bigl(\cp_n\cup \cp_n^{(k)}\bigr)$, one naturally suppose that $L_n$ is a subset of $N_n^{(k)\natural}$ such that $L_n\cap \part N_n^{(k)\natural}=\part_\cs N_n^{(k)\natural}$.
Note that the inclusion of $L_n$ to $N_n^{(k)\natural}$ is $\fd$-injective.
Let $p_n^{(k)}:\tilde N_n^{(k)}\to N_n^{(k)}$ be the locally isometric covering associated with $\fd(L_n)\subset \fd(N_n^{(k)\natural})$ and $\rho_n^{(k)}:\fd(L_n)\to \mathrm{PSL}_2(\cc)$ the holonomy of $\tilde N_n^{(k)}$.
In particular, $\tilde N_n^{(k)}$ is extended to a Kleinian manifold $\tilde N_n^{(k)\natural}$ which can be identified with 
$\mathbf{H}^3\cup \O(\G_n^{(k)})/\G_n^{(k)}$, where $\G_n^{(k)}=\mathrm{Image}(\rho_n^{(k)})$.
Let $\tilde L_n^{(k)}$ be a proper core of $\tilde N_n^{(k)\natural}$ such that 
$p_n^{(k)\natural}|\tilde L_n^{(k)}$ is a proper embedding to $L_n$ such that 
$p_n^{(k)\natural}(\tilde L_n^{(k)})$ is a proper deformation retract of $L_n$.

\begin{lemma}\label{s_convergence}
The $k$-sequence $\{\rho_n^{(k)}\}$ subconverges strongly to a Kleinian representation $\rho_n$ as $k\sto \infty$.
\end{lemma}
\begin{proof}
From our construction, the end of $\tilde N_n^{(k)}$ corresponding to $F_j$ has the geodesic loop $\tilde l_j^{(k)}$ with $p_n^{(k)}(\tilde l_j^{(k)})=\ol l_j^{(k)}$.
From the fact that $\lim_{k\sto \infty}\mathrm{length}(\tilde l_j^{(k)})=0$ together with arguments in Ohshika 
\cite{oh0},  the sequence $\{\rho_n^{(k)}\}$ subconverges algebraically to a Kleinian representation $\rho_n$ 
as $k\sto \infty$ such that the end of $\tilde N_n=\mathbf{H}^3/\G_n$ corresponding to $F_j$ is geometrically infinite and the its ending lamination coincides with the support of $\mu_j$, where $\rho_n(\fd(L_n))=\G_n$.
Moreover, $\part \tilde N_n^\natural =\part_\cs \tilde N_n^\natural$ is conformally equivalent 
to $\part_\cs \tilde N_n^{(k)\natural}$ for all sufficiently large $k$.
Then one can show that the $k$-sequence $\{\rho_n^{(k)}\}$ subconverges strongly to $\rho_n$ by applying Thurston's Covering Theorem \cite[Theorem 9.2.2]{tl}.
\end{proof}

Since $L_n$ is a $\fd$-injectively embedded subset of $L_m$ for any $m,n\in \nn$ with $n\leq m$, one can suppose that 
$\fd(L_n)$ is a subgroup of $\fd(L_m)$.
The subgroup $\rho_m(\fd(L_n))$ of $\G_m=\rho_m(\fd(L_m))$ is denoted by $\G_{m;n}$, and the 
restriction $\rho_m |\fd(L_n)$ by $\rho_{m;n}$.
Let $p_{m;n}:\tilde N_{m;n}=\mathbf{H}^3/\G_{m;n}\to \tilde N_m$ be the locally isometric 
covering associated with $\G_{m;n}\subset \G_m$.

\begin{lemma}\label{a_convergence2}
For each $n\in \nn$, the $m$-sequence $\{\rho_{m;n}\}_{m\geq n}$ subconverges algebraically 
to a Kleinian representation $\rho_{\infty;n}$ as $m\sto \infty$.
\end{lemma}

\begin{proof}
For any non-peripheral element $F_j$ of $\cz_n$, let $f_{m;j}:F_j(\sigma_{m;j})\to \tilde N_{m;n}$ be an efficient pleated surface facing to the end of $\tilde N_{m;n}$ corresponding to $F_j$.
Here we will show that the closure of the $m$-sequence $\{\sigma_{m;j}\}$ in $\mathcal{ML}_0(F_j)$ is compact.
For the proof, suppose contrarily that $\{\sigma_{m;j}\}$ subconverges to $[\mu']\in \partial
\mathrm{Teich}(F_j)= \mathcal{PL}_0(F_j)$ as $m\sto \infty$ 
and introduce a contradiction.

Fix a component $\mu$ of $\mu'$.
When $F_j\subset \Sg_v$, let $A^-(\mu)=\mu_{[u,v]}$, $A^+(\mu)=\mu_{[v,w]}$ be subspaces of $\Sg\times I$ such that 
(i) $\mu_{(u,v)}$, $\mu_{(v,w)}$ are subsets of $L$ and (ii) $A^\pm(\mu)$ are 
maximal among all $\mu_{[\cdot,v]}$, $\mu_{[v,\cdot]}$ satisfying (i).
Set $\mu^-=A^-(\mu)\cap \Sg_u$ and $\mu^+=A^+(\mu)\cap \Sg_w$.
One can suppose that $\mu^+$ is neither contained in an element of $\cp_\cx$ nor equal to the fixed lamination $\mu_t$ of a non-peripheral front $F_t$ of $\cx$, otherwise we may use $\mu^-$ 
instead of $\mu^+$.
It suffices to consider the following three cases.

{\bf Case 1.}
$\mu^+$ meets a component $P$ of $\mathcal{P}_\cx$ nontrivially.

From our assumption as above, $\mu^+$ meets $P$ transversely.
By Lemma \ref{c_to_b}, $A^+(\mu)$ passes through only finitely many proper vertical surfaces $H_l$ $(l=0,1,\dots,2a)$ in $L$ with 
$H_0=F_j$ such that $\mathrm{pr}_{\mathrm{hz}}(H_{2i-1})$ is a component of 
$\mathrm{pr}_{\mathrm{hz}}(H_{2i-2})\cap \mathrm{pr}_{\mathrm{hz}}(H_{2i})$ for $i=0,1,\dots,a$, see Fig.\ \ref{f_5_1}.
\begin{figure}[hbtp]
\centering
\scalebox{0.5}{\includegraphics[clip]{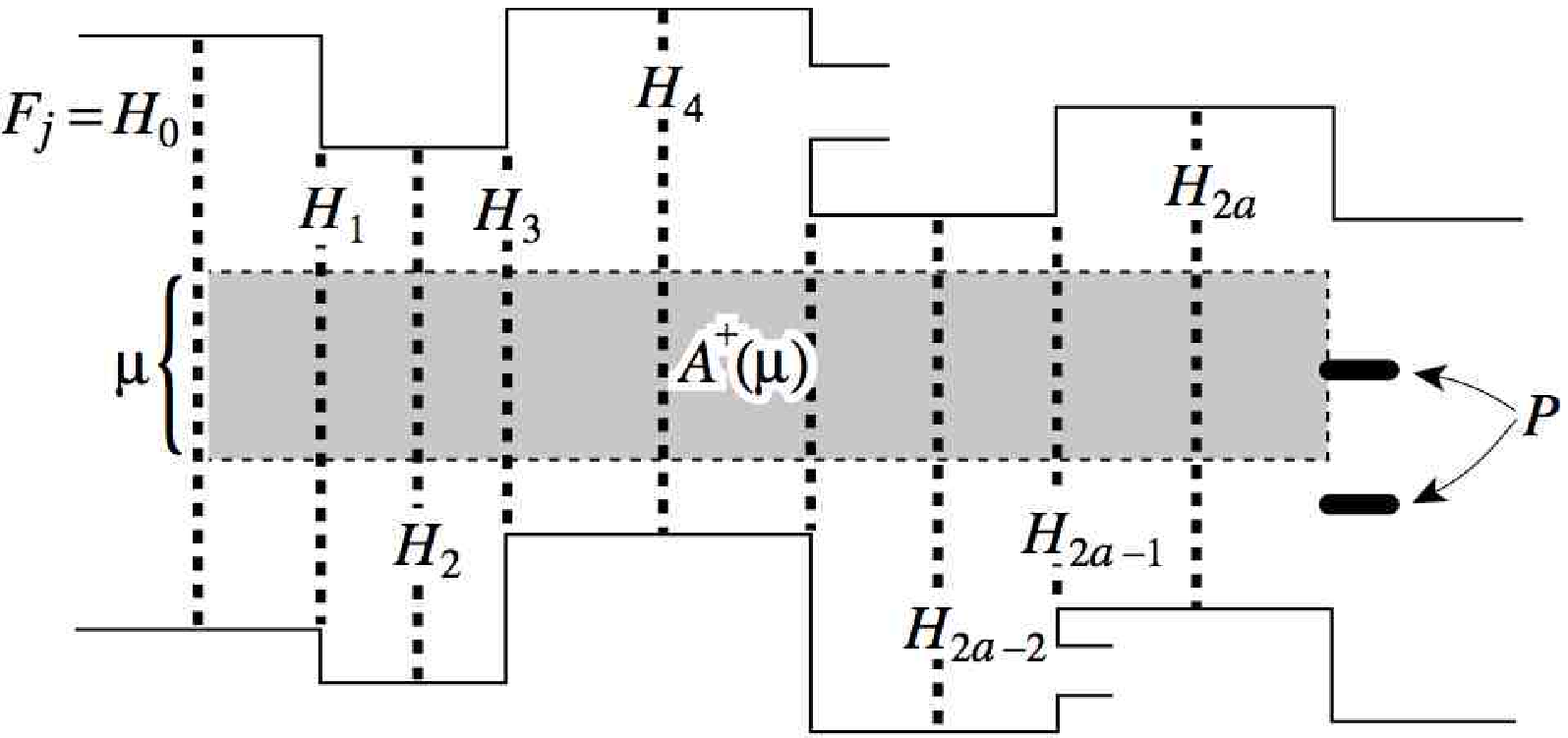}}
\caption{}
\label{f_5_1}
\end{figure}
For all $m$ sufficiently larger than $n$, $L_m$ contains $H_0\cup\dots\cup H_{2a}$.
Let $g_{m;i}:H_i(\tau_{m;i})\to \tilde N_{m}$ $(i=0,1,\dots,2a)$ be an efficient pleated surface corresponding to $H_i\subset L_m$ with 
$g_{m;0}=p_{m;n}\circ f_{m;j}$ and $\tau_{m;0}=\sigma_{m;j}$.
Since $\mathrm{length}_{\tilde N_{m}}(P)=\mathrm{length}_{\tilde N_{m}}(\fr(H_{2a-1}))=0$, by (\ref{length}), both $\mathrm{length}_{\tau_{m;2a}}(P)$ and $\mathrm{length}_{\tau_{m;2a}}(c_{2a-1})$ are uniformly bounded by constants independent of $m$, where $c_{2a-1}$ is the disjoint union of g-loops in $H_{2a}$ parallel to $\fr(H_{2a-1})$ in $L$.
By Lemma \ref{l5_1}, there exists a measured lamination $\lambda_{2a}$ in $H_{2a}$ with 
$\mathrm{pr}_{\mathrm{hz}}(\lambda_{2a})\subset \mathrm{pr}_{\mathrm{hz}}(H_{2a-1})$ and such that 
$\mathrm{length}_{\tau_{m;2a}}(\lambda_{2a})$ is uniformly bounded and $i(\mu,\lambda_{2a})\neq 0$.
Since $\mathrm{pr}_{\mathrm{hz}}(H_{2a-2})\supset \mathrm{pr}_{\mathrm{hz}}(H_{2a-1})$, $H_{2a-2}$ 
has a measured lamination $\lambda_{2a-2}$ with $\mathrm{pr}_{\mathrm{hz}}(\lambda_{2a-2})=
\mathrm{pr}_{\mathrm{hz}}(\lambda_{2a})$.
Again by (\ref{length}), $\mathrm{length}_{\tau_{m;2a-2}}(\lambda_{2a-2})$ is uniformly 
bounded. 
By repeating the same argument $a$ times, we have a measured laminations $\lambda_0$ in $H_0=F_j$ such that $\mathrm{length}_{\sigma_{m;j}}(\lambda_0)$ is uniformly bounded by a constant independent of $m$ and $i(\mu,\lambda_0)\neq 0$.
This contradicts that $\{\sigma_{m;j}\}$ subconverges to $[\mu]$ as $m\sto \infty$, for example 
see \cite[Theorem 2.2]{tii}.

{\bf Case 2.}
$\mu^+$ is contained in a non-peripheral front $F_t$ of $\cx$ nontrivially.

Since the support of $\mu_t$ is the ending lamination of the end of 
$\tilde N_{m}$ corresponding to $F_t$ for all $m$, $\mathrm{length}_{\tilde N_{m}}(\mu_t)=0$.
Thus one can introduce a contradiction by an argument quite 
similar to that in Case 1.

{\bf Case 3.}
$\mu^+$ is contained in $\Sg\times \{-1,1\}$.

Then $\mu^+$ is a measured lamination of some peripheral font $S$ of $\cx$.
Since the conformal structure on the component of $\partial_{\mathcal{S}} \tilde N_{m}^\natural$ 
corresponding to $S$ is fixed, the length of any geodesic loop in $\tilde N_{m}$ parallel to a simple loop in $S$ is uniformly bounded.
Thus, also in this case, one can obtain a contradiction.

The results in Cases 1-3 imply that, for each $j$, the closure of the $m$-sequence $\{\sigma_{m;j}\}$ in $\mathcal{ML}_0(F_j)$ is compact.
By Thurston's Relative Boundedness Theorem \cite[Theorem 3.1]{tiii}, the $m$-sequence $\{\rho_{m;n}\}_{m\geq n}$ subconverges algebraically to a Kleinian representation $\rho_{\infty;n}$.
\end{proof}

Suppose that $\{\xi_n:\Lambda\to \mathrm{PSL}_2(\cc)\}$ is a sequence of Kleinian representations 
converging algebraically to $\xi$.
We say that $\xi(\lambda)$ with $\lambda\in \Lambda$ is an \emph{accidental parabolic element} 
with respect to $\{\xi_n\}$ if $\xi(\lambda)$ itself is parabolic, but $\xi_n(\lambda)$ is not so for all sufficiently 
large $n\in\nn$.

\begin{lemma}\label{s_convergence2}
The convergence $\{\rho_{m;n}\}_{m\geq n}\sto \rho_{\infty;n}$ as $m\sto \infty$ is strong.
\end{lemma}
\begin{proof}
If not, by Anderson and Canary \cite[Corollary 3.3]{ac2}, there would exist an element $\gamma\in \fd(L_n)$ such that $\rho_{\infty;n}(\gamma)$ is an accidental parabolic element.
It is well known that the element $\gamma$ is represented by a g-loop $l_n$ in some element $F_j$ of $\mathcal{Z}_n$.
Since any non-peripheral element of $\mathcal{Z}_n\setminus \mathcal{J}_n$ corresponds to a geometrically infinite end of $\tilde N_{\infty,n}$ with an ending lamination without compact 
leaves, 
$F_j$ either is peripheral or belongs to the joint set $\mathcal{J}_n$.
Since $\rho_{\bar m;n}=\rho_{\bar m;m}|\fd(L_n)$ for any $\bar m\geq m\geq n+1$, $\rho_{\infty;n}=\rho_{\infty;m}|\fd(L_n)$.
Hence, for any $m\geq n$, $\tilde N_{\infty,m}$ has the accidental parabolic $\mathbf{Z}$-cusp 
corresponding to a g-loop $l_m$ which is parallel to $l_n$ in $L_m$ and either peripheral or contained in some element of $\mathcal{J}_m$.
Since our block complex satisfies the finiteness property given in Lemma \ref{c_to_b}, by 
repeating this argument finitely many times, one can show that 
$l_n$ is parallel in $L$ to a g-loop $l_\infty$ in a component $S$ of $\cs$.
Since the marked conformal structure on the component $S_{m;n}$ of $\part \tilde N_{m;n}^\natural$ corresponding to 
$S$ is fixed, $\part \tilde N_{\infty;n}^\natural$ also has the component $S_{\infty;n}$ with the same marked conformal 
structure.
In particular, $N_{\infty;n}^\natural$ has no parabolic cusps adjacent to $S_{\infty;n}$.
This contradiction implies our desired strong convergence.
\end{proof}

\subsection{Proof of Theorem \ref{main_b}}\label{ss_5_4}

Since $\rho_{\infty;n}=\rho_{\infty;n+1}|\fd(L_n)$, $\G_{\infty;n}=\rho_{\infty;n}(\fd(L_n))$ is a subgroup of $\G_{\infty;n+1}=\rho_{\infty;n+1}(\fd(L_{n+1}))$ for any $n\in \nn$.
This implies that the union $G=\bigcup_{n\geq 1}\G_{\infty;n}$ is a geometric limits of the $n$-sequence $\{\G_{\infty;n}\}$.
If necessary passing to a subsequence, we may assume that $\{\G_m\}$ converges geometrically to a Kleinian group $G'$ as $m\sto \infty$.
Since the $m$-sequence $\{\G_{m;n}\}$ converges algebraically to $\G_{\infty;n}$ and $\G_{m;n}\subset \G_m$, 
$\G_{\infty;n}$ is a subgroup of $G'$ for any $n\geq 1$ and hence $G$ is also a subgroup of $G'$.

The proof of Theorem \ref{main_b} is completed by the following three propositions.

\begin{prop}\label{main_b_1}
$M^\natural=\mathbf{H}^3\cup \O(G)/G$ is homeomorphic to $L=\Sg\times I\setminus \cx$.
\end{prop}
\begin{proof}
We use standard arguments of the least area surface theory in \cite{fhs} again.
Since the $m$-sequence $\{\rho_{m;n}\}$ converges strongly to $\rho_{\infty,n}$ by Lemma \ref{s_convergence2}, 
$\tilde N_{\infty;n}^\natural$ has a proper core $\tilde L_{\infty;n}$ such that $\tilde L_{\infty;n}\cap \tilde N_{\infty;n,\pthick}^\natural$ 
is mapped via a marking-preserving bi-Lipschitz map to a submanifold of $\tilde N_{m;n,\pthick}^\natural$ properly homotopic to 
$\tilde L_{m;n}\cap \tilde N_{m;n,\pthick}^\natural$, where $\tilde L_{m;n}$ is a proper core of $\tilde N_{m;n}^\natural$. 
For any positive integers $n,\bar n$ with $n\leq \bar n$, let $q_{\bar n; n}:\tilde N_{\infty;n}\to \tilde N_{\infty;\bar n}$ 
be the covering associated with $\G_{\infty;n}\subset \G_{\infty;\bar n}$ and $q_{\infty;n}:\tilde N_{\infty;n}\to M$ 
the covering associated with $\G_{\infty;n}\subset G$.
Modify the hyperbolic metric on $M_{\mathrm{p:thick}}$ in a small collar neighborhood of $\part M_{\mathrm{p:thick}}$ so that $\part M_{\mathrm{p:thick}}$ has non-negative mean curvature in $M_{\mathrm{p:thick}}$.
Suppose that $W_n=q_{\infty;n}^{-1}(M_{\mathrm{p:thick}})$ has the Riemannian metric induced from that on 
$M_{\mathrm{p:thick}}$.
Since $W_n$ is a deformation retract of $\tilde N_{n;\infty}$, we may assume that $W_n\cap \tilde L_{n;\infty}$ is a compact 
core of $\tilde L_{n;\infty}$ and $\tilde L_{n;\infty}\setminus W_n$ consists of parabolic cusps.
By \cite{fhs}, we may also assume that, for any non-peripheral component $F_{n;j}$ of $\fr(\tilde L_{n;\infty})$, $F_{n;j}\cap W_n$ is a least 
area surface in its proper homotopy class.
For any $\bar n\geq n$ and $F_{\bar n,j}=q_{\bar n;n}(F_{n;j})$, $F_{\bar n;j}\cap W_{\bar n}$ are also least area surfaces in $W_{\bar n}$ properly homotopic to mutually disjoint embedded surfaces in $\tilde L_{\bar n;\infty}\cap W_{\bar n}$.
Again by \cite{fhs}, $F_{\bar n;j}$ are also mutually disjoint embedded surfaces in $\tilde N_{\bar n;\infty}$ and hence 
$F_{\infty;j}=q_{\infty;n}(F_{n;j})$ are so in $M$.
Then the family $\{F_{\infty;j}\}_{j}$ determines a block decomposition on $M^\natural$ such that the support $M^\natural$ 
of the decomposition is homeomorphic to $L$.
\end{proof}

\begin{prop}\label{main_b_2}
$G=G'$.
\end{prop}
\begin{proof}
Let $p: M\to M'$ be the covering associated with $G\subset G'$, where $M'=\mathbf{H}^3/G'$.
By Proposition \ref{main_b_1}, $\part M^\natural=\part_\cs M^\natural$.
Since $\part_\cs \tilde N_m^\natural$ is the union of components of $\part \tilde N_m^\natural$ 
which exist permanently for all $m$, we have $\part {M'}^\natural=\part_\cs {M'}^\natural$.
Hence, if $\cs$ is nonempty, then $p$ extends to a homeomorphism from $\part M^\natural$ to $\part {M'}^\natural$.
It follows that $p$ is also homeomorphism.
Thus it suffices to consider the case when $\cs=\eset$, or equivalently both $M,M'$ have no 
geometrically finite ends.

Take a parabolic cusp $P_0$ of $M$.
For any sufficiently large $n$, $\tilde N_{\infty;n}$ has a unique parabolic cusp $P_1$ covering $P_0$.
Since $\G_{\infty;n}$ has no accidental parabolic elements with respect to the strong convergence of $\{\rho_{m;n}\}$, 
for any sufficiently $m$, $\tilde N_{m;n}$ has a unique cusp $P_2$ corresponding to $P_1$ which in turn maps to a parabolic 
cusp of $P_3$ of $\tilde N_m$.
Since $\{\G_m\}$ converges geometrically to $G'$, $M'$ has a unique cusp $P_4$ corresponding to $P_3$.
Since the inclusion $\G_{m;n}\subset \G_m$ is derived from $L_n\subset L_m$, the correspondence $P_1\mapsto P_2$ are injective.
Moreover, it follows from the fact that the other correspondences are defined via bi-Lipschitz maps, 
that distinct parabolic cusps $P_0$ of $M$ correspond to distinct parabolic cusps $P_4$ of $M'$. 
It follows that $p:M\to M'$ defines an embedding from the union of parabolic cusps of $M$ into that of $M'$.
So, if $p^{-1}(P_4)$ had a component $B$ other than $P_0$, then $B$ would be a horoball.
Since $M$ does not have geometrically finite ends, there would exist a pleated surface $F$ in $M$ passing through a point of $\Int B$ arbitrarily 
far from $\part B$ by the ubiquity of pleated surfaces in convex cores.
This implies that $F$ would contain a hyperbolic disk of arbitrarily large diameter, a contradiction.
Thus the restriction of $p$ on $p^{-1}(P_4)$ is a homeomorphism onto $P_4$, and hence the covering $p$ is also a homeomorphism.
\end{proof}

\begin{prop}\label{main_b_3}
$G'$ is a geometric limit of algebraically convergent quasi-Fuchsian groups.
\end{prop}
\begin{proof}
The group $G'$ is geometrically approximated by $\G_m$ and hence by $\G_m^{(k)}$ for all sufficiently large $m$ 
and $k$.
Take an integer $R>0$ arbitrarily.
Any element $T_l$ of $\mathcal{T}_m$ is realized as a totally geodesic surface in $M'$.
For all sufficiently large $m$, all such totally geodesic surfaces do not intersect $\cb_R(x_\infty',M')$.
Pleated surfaces $F_j^{(k)}$ realizing the g-loop $l_j^{(k)}$ for $j=1,\dots,s$ and totally geodesic surfaces $T_l^{(k)}$ in $N_m^{(k)}$ corresponding to $T_l\in \mathcal{T}_m$ eventually leave $\cb_R(x_{m,k},N_m^{(k)})$ as $k,m\sto \infty$.
Since the component of $N_m^{(k)}\setminus \bigl(\bigcup F_j^{(k)}\bigr)\cup\bigl(\bigcup T_l^{(k)})$ containing $x_{m,k}$ is lifted to $\tilde N_j^{(k)}$ (see Fig.\ \ref{add_1}), 
the restriction of $p_m^{(k)}$ to $\cb_R(\tilde x_{m,k},\tilde N_m^{(k)})$ is an isomorphism onto $\cb_R(x_{m,k},N_m^{(k)})$.
This implies that $\G_m^{(k)}$ is geometrically approximated by a geometrically finite Kleinian group $\Lambda_m^{(k)}$ such that $N_m^{(k)}$ is isometric to $\mathbf{H}^3/\Lambda_m^{(k)}$.
Applying the Maskit Second Combination Theorem \cite[Theorem E.5]{mas}, one can obtain a new hyperbolic 3-manifold 
$\widehat N_m^{(k)}=\mathbf{H}^3/\widehat \Lambda_m^{(k)}$ which replaces each $\mathbf{Z}$-cusp of $N_m^{(k)}$ by a $\mathbf{Z}\times \mathbf{Z}$-cusp 
without changing the part $\cb_R(x_m^{(k)},N_m^{(k)})$ so that the new manifold is homeomorphic to
$\Sg\times I$ minus finitely many g-loops in $\Int (\Sg\times I)$.
By performing hyperbolic Dehn surgeries along the parabolic cusps of $\widehat N_m^{(k)}$ of type $(1,u_k)$ with sufficiently large 
$u_k$, we know that $\widehat \Lambda_m^{(k)}$ and hence $\Lambda_m^{(k)}$ are geometrically approximated by a quasi-Fuchsian group 
$\L_R$ admitting a representation $\zeta_R:\Pi\to \mathrm{PSL}_2(\cc)$.
This shows that $\{\L_R\}$ converges geometrically to $G$ as $R\sto \infty$.
From our construction, $\{\zeta_R\}$ converges algebraically to the restriction 
$\rho_{\infty;n}|\Pi$, which is independent of $n$.
This completes the proof of Proposition \ref{main_b_3} and hence that of Theorem \ref{main_b}.
\end{proof}


\begin{thebibliography}{99}
\bibitem{tii} W. Thurston, Hyperbolic structures on $3$-manifolds, II: Surface groups and $3$-manifolds which fiber over the circle, E-print: {\tt math.GT/9801045}.
\end{thebibliography}
\end{document}